\documentclass[a4paper]{article}

\def\FLOPs{flop/s} 
\newcommand{\tlapack}{{$\langle$T$\rangle$LAPACK}}
\newcommand{\cpp}{{C\texttt{++}}}

\usepackage{utfsym}
\definecolor{emerald}{rgb}{0.31, 0.78, 0.47}
\newcommand{\yes}{\textcolor{emerald}{\usym{2714}}}
\newcommand{\no}{\textcolor{red}{\usym{2718}}}
\newcommand{\maybe}{\textcolor{orange}{\usym{1F311}}}
\newcommand{\ignore}[1]{}
\usepackage{array}

\usepackage[T1]{fontenc}
\usepackage[utf8]{inputenc}
\usepackage{lmodern}
\usepackage{microtype}
\usepackage{amsmath,amssymb}
\usepackage{enumitem}
\usepackage{doi}
\usepackage{csquotes}
\usepackage{verbatim}
\usepackage{dirtytalk}
\usepackage{booktabs}
\PassOptionsToPackage{hyphens}{url}
\usepackage{hyperref}
\usepackage{geometry}
\hypersetup{
  colorlinks=true,
  breaklinks=true,
  linkcolor=blue,
  urlcolor=blue,
  citecolor=blue
}
\usepackage{cleveref, soul,todonotes}

\usepackage[utf8]{inputenc}
\usepackage[acronym]{glossaries}

\usepackage{tikz}
\usepackage[dvipsnames]{xcolor}
\usetikzlibrary{positioning, arrows.meta, calc, decorations.pathreplacing}

\definecolor{hwColor}{HTML}{2E4053}
\definecolor{isaColor}{HTML}{34495E}
\definecolor{rtlColor}{HTML}{5D6D7E}
\definecolor{drvColor}{HTML}{2874A6}
\definecolor{runtimeColor}{HTML}{1E8449}
\definecolor{fwkColor}{HTML}{D4AC0D}
\definecolor{appColor}{HTML}{A04000}
\definecolor{ifaceColor}{HTML}{7D3C98}
\definecolor{mpColor}{HTML}{C0392B}
\definecolor{toolColor}{HTML}{2E86C1} %

\makeglossaries
\newacronym{AGQ}{AGQ}{AMD-GraphCore-Qualcomm Floating Point Format}
\newacronym{AI}{AI}{Artificial Intelligence}
\newacronym{ALU}{ALU}{Arithmetic-Logic Unit}
\newacronym{ALUs}{ALUs}{Arithmetic-Logic Units}
\newacronym{AMD}{AMD}{Advanced Micro Devices}
\newacronym{AMR}{AMR}{Adaptive Mesh Refinement}
\newacronym{API}{API}{Application Program Interface}
\newacronym{ASIC}{ASIC}{Application-Specific Integrated Circuit}
\newacronym{BF16}{BF16}{Brain Float16 Floating Point Format}
\newacronym{BFP}{BFP}{Block Floating Point}
\newacronym{BLAS}{BLAS}{Basic Linear Algebra Subprograms}
\newacronym{BLR}{BLR}{Block Low Rank}
\newacronym{CESTAC}{CESTAC}{Control and Stochastic Estimation of Rounding Errors}
\newacronym{CG}{CG}{Conjugate Gradient}
\newacronym{CMG}{CMG}{Compact Multigrid}
\newacronym{CPU}{CPU}{Central Processing Unit}
\newacronym{DGEMM}{DGEMM}{Double-Precision GEneral Matrix Multiply}
\newacronym{DRAM}{DRAM}{Dynamic Random Access Memory}
\newacronym{DSA}{DSA}{Discrete Stochastic Arithmetic}
\newacronym{EPI}{EPI}{European Processor Initiative}
\newacronym{FMA}{FMA}{Fused Multiply-Add}
\newacronym{FPGA}{FPGA}{Field Programmable Gate Array}
\newacronym{FPU}{FPU}{Floating Point Unit}
\newacronym{FPUs}{FPUs}{Floating Point Units}
\newacronym{FP8}{FP8}{8-bit floating point} 
\newacronym{FP16}{FP16}{IEEE~754 \texttt{binary16}} 
\newacronym{FP32}{FP32}{IEEE~754 \texttt{binary32}} 
\newacronym{FP64}{FP64}{IEEE~754 \texttt{binary64}} 
\newacronym{GEMM}{GEMM}{GEneral Matrix Multiply}
\newacronym{GMRES}{GMRES}{the Generalized Minimal Residual Method}
\newacronym{GPU}{GPU}{Graphics Processing Unit}
\newacronym{HiF8}{HiF8}{Huawei HiFloat8 Floating Point Format}
\newacronym{HODLR}{HODLR}{Hierarchically Off-Diagonal Low Rank}
\newacronym{HPC}{HPC}{High-performance computing}
\newacronym{HSS}{HSS}{Hierarchical Semi-Separable}
\newacronym{IR}{IR}{Iterative Refinement}
\newacronym{IEEE}{IEEE}{Institute of Electrical and Electronics Engineers}
\newacronym{INT8}{INT8}{8-bit signed integer format using two's complement}
\newacronym{INT32}{INT32}{32-bit signed integer format using two's complement} 
\newacronym{ISA}{ISA}{Instruction Set Architecture}
\newacronym{LNSM}{LNSM}{NVIDIA LSN-Madam Floating Point Format}
\newacronym{MCA}{MCA}{Monte-Carlo Arithmetic}
\newacronym{MG}{MG}{Multigrid}
\newacronym{MMA}{MMA}{Matrix Multiply-Accumulate}
\newacronym{MPI}{MPI}{Message Passing Interface}
\newacronym{MPFR}{MPFR}{GNU C library for multiple-precision floating-point computations with correct rounding}
\newacronym{NaN}{NaN}{Not-a-Number}
\newacronym{OCP}{OCP}{Open Compute Project}
\newacronym{MLP}{MLP}{Multilayer Perceptron}
\newacronym{MXP}{MX}{Open Compute Project Floating Point Format}
\newacronym{PDE}{PDE}{Partial Differential Equation}
\newacronym{Posit}{Posit}{Posit Floating Point Format}
\newacronym{P3109}{P3109}{IEEE P3109 Working Group on Binary Floating Point Formats for Machine Learning}
\newacronym{QoI}{QoI}{Quantity of Interest}
\newacronym{RN}{RN}{Round to Nearest}
\newacronym{RU}{RU}{Round Up}
\newacronym{RD}{RD}{Round Down}
\newacronym{RDS}{RDS}{Rounding Deferred Summation}
\newacronym{RZ}{RZ}{Round toward Zero}
\newacronym{RTO}{RTO}{Round to Odd}
\newacronym{SEM}{SEM}{Spectral Element Methods}
\newacronym{SIMD}{SIMD}{Single Instruction Multiple Data}
\newacronym{SMT}{SMT}{Satisfiability Modulo Theories}
\newacronym{SpMV}{SpMV}{Sparse Matrix-Vector Multiplication}
\newacronym{SR}{SR}{Stochastic Rounding}
\newacronym{SVD}{SVD}{Singular Value Decomposition}
\newacronym{Takum}{Takum}{Takum Floating Point Format}
\newacronym{TF32}{TF32}{Tensor Flow32 Floating Point Format}
\newacronym{TSL}{TSL}{Tesla Dojo Floating Point Format}
\newacronym{UQ}{UQ}{Uncertainty Quantification}
\newacronym{VaV}{V\&V}{Verification \& Validation}

\newacronym{754}{754}{IEEE 754 Floating Point Standard}

\begin{document}
    
\title{Mixed-Precision Computing for Scientific Discovery: Formats, Co-Design, and Responsible Approximation}

\author{%
\small
  \begin{tabular}{c}
    Emmanuel Agullo (INRIA - Bordeaux, FR), Hartwig Anzt (TU München - Heilbronn, DE),
    Daniel Bauer\\ (Universität Erlangen-Nürnberg, DE), David Bindel (Cornell University - Ithaca, US), 
    Alfredo Buttari (CNRS\\ IRIT - Toulouse, FR), Alexandru Calotoiu (ETH Zürich, CH), 
    Erin Claire Carson
    (Charles University,  CZ),\\
    Pasqua D'Ambra (CNR IAC - Naples, IT), Ieva Daužickait\.{e} (Cerfacs - Toulouse, FR), 
    James W. Demmel\\ (University of California - Berkeley, US), 
    Jack Dongarra (University of Tennessee - Knoxville, US),
    Iain Duff\\ (STFC Rutherford Appleton Laboratory, GB \& Cerfacs, FR), 
    Massimiliano Fasi (University of Leeds, GB),\\
    Dominik Göddeke (Universität Stuttgart, DE),    Stef Graillat (Sorbonne University - Paris, FR),
    Laslo Hunhold\\ (Openchip - Barcelona, ES),
    Roman Iakymchuk (Umeå University, SE \& Uppsala University, SE),
    Fabienne\\ Jézéquel (Sorbonne University - Paris, FR),
    Nils Kohl (LMU München, DE),
    Harald Köstler (Universität\\ Erlangen-Nürnberg, DE),
    Jakub Kružík (Technical University of Ostrava \& Czech Academy of Sciences,\\ Institute of Geonics, CZ),
    Julien Langou (University of Colorado Denver - Denver, US),
    Xiaoye Sherry Li\\ (Lawrence Berkeley National Laboratory, US),
    Hatem Ltaief (KAUST - Thuwal, SA),
    Piotr Luszczek (MIT\\ Lincoln Laboratory, US),
    Yuxin Ma (Charles University, CZ),
    Theo Mary (Sorbonne Université, CNRS,\\ LIP6 - Paris, FR),
    Mantas Mikaitis (University of Leeds, GB),
    Hiroyuki Ootomo (NVIDIA - Tokyo, JP),\\
    Daniel Osei-Kuffuor (LLNL - Livermore, US),
    Enrique S. Quintana-Ortí (Universitat Politècnica de València,\\ ES),
    Ulrich Rüde (Universität Erlangen-Nürnberg, DE),
    Jennifer Scott (University of Reading, GB \& Rutherford\\ Appleton Laboratory, GB),
    John Shalf (Lawrence Berkeley National Laboratory, US),
    Linda Stals (Australian\\ National University - Canberra, AU),
    Rasmus Tamstorf (Intel Corporation, US),
    Stefan Turek (TU Dortmund,\\ DE),
    Petr Vacek (IFP Energies nouvelles, FR),
    Bastien Vieublé (Chinese Academy of Sciences - Beijing, CN),\\
    Rio Yokota (Institute of Science Tokyo, JP)
  \end{tabular}%
}

\date{}

\maketitle

\begin{abstract}
Reduced and mixed precision have moved from a niche optimization to a central design axis in scientific computing and engineering, driven by energy constraints, heterogeneous accelerators, and the convergence of simulation and machine learning. This paper organizes the landscape around seven coupled themes---number formats, floating-point emulation, emerging architectures, hardware/software co-design, relation to other approximations, software design, and precision as a multilevel resource ---and, for each theme, synthesizes the state of the art, future directions, and open questions. We emphasize \emph{energy per trusted solution} as the core objective, and we frame \say{recklessly responsible} computing as a pragmatic doctrine: exploit low precision aggressively, but with systematic detection, escalation, and certification pathways.
\footnotetext{This work is dedicated to the memory of Nicholas J. Higham in appreciation of his lasting contributions to the area of numerical computing.}

\end{abstract}

\section{Introduction: why precision is now a first-class design variable}

For decades, computational science and engineering largely assumed 
 \acrfull{FP32}  and \acrfull{FP64}
as the numerical baseline, turning to higher precision only in exceptional situations. That assumption no longer holds. Contemporary accelerators provide a widening spectrum of low-precision formats and  --just as importantly-- a variety of accumulation and rounding behaviors with sharply different performance, energy characteristics, and numerical properties. 
At the same time, modern workflows increasingly layer multiple sources of error: modeling and discretization error, iterative and sampling error, data uncertainty, and roundoff. 
These effects often coexist within a single end-to-end pipeline, including data-driven components.

As a consequence, the central objective changes. The long-held goal of \say{\FLOPs\ at a fixed precision} has been replaced by \emph{minimizing time and energy subject to a trust target}. That trust target may be expressed as an accuracy tolerance for a quantity of interest, a validated uncertainty interval, a reproducibility requirement, or a failure-detection criterion.

This paper offers a practical way to organize the mixed-precision 
landscape in seven themes: number formats, floating-point emulation, emerging architectures, hardware/software co-design, relation to other approximations, software design, and precision as a multilevel resource.
For each topic area, we distinguish:
\begin{itemize}[leftmargin=1.5em]
  \item \textbf{State of the art:} techniques that are dependable in current systems.
  \item \textbf{Future directions:} plausible advances over the next 3--10 years.
  \item \textbf{Open questions:} barriers to standardization, portability, and broader adoption.
\end{itemize}

Across topics, four goals recur:
\begin{itemize}[leftmargin=1.5em]
  \item \textbf{Trust:} \acrfull{VaV}, \acrfull{UQ}, reproducibility, and robust failure detection.
  \item \textbf{Portability:} comparable behavior across hardware, compilers, and libraries.
  \item \textbf{Composability:} precision decisions that remain valid under other numerical approximations, multi-physics coupling, and coupled \acrshort{AI}/simulation pipelines.
  \item \textbf{Efficiency:} evaluated as time or energy per \emph{trusted} solution, not peak throughput.
\end{itemize}

\section{Floating-point number formats}\label{sec:numberformats}

\label{sec:NumberFormats}

Floating-point number formats are quantization schemes for the infinite set of real numbers, defining a predefined, finite set of representatives that are used to approximate the true numbers.
Since no finite set of reals can be closed under the arithmetic operations, this leads to the questions of rounding, overflow and underflow.
Currently, many new such schemes are under consideration with the goal of permitting better tradeoffs between precision and cost---typically reducing cost by sacrificing precision.
The resulting landscape of current reduced-precision formats is both rich and messy. 
On the one hand, commercially deployed systems now expose multiple low-precision formats and specialized execution units that deliver large gains in throughput and energy efficiency. 
On the other hand, the formats themselves---their exception semantics, rounding behavior, saturation options (how overflow is handled), and conversion rules---vary across vendors and are often only partially specified in public. 
As a result, choosing a format for scientific work is rarely a matter of \say{bits-per-value} alone; it is a question of \emph{what behavior is guaranteed}, \emph{what can be implemented portably}, and \emph{what can be analyzed and certified}.

\subsection{State of the art}

To make this discussion actionable, we evaluate formats against a set of desiderata that reflect the needs of numerical computing beyond machine learning. These criteria span practical deployment questions (hardware availability, patent encumbrance, portability), semantic questions (uniqueness of representations, well-defined handling of exceptional values, saturation and rounding modes), and mathematical questions (whether existing numerical analysis applies, whether higher precision can be composed/emulated reliably, and whether basic algebraic properties hold). The goal is not to crown a single \say{best} format, but to clarify the trade-offs that determine whether a format can serve as a \emph{portable primitive} for mixed-precision scientific software.

We then survey a representative set of formats that appear in today’s ecosystem: \acrshort{IEEE}~754 and its derivatives (\acrshort{BF16}, \acrshort{TF32}), emerging \acrshort{AI}-oriented proposals and consortia definitions (e.g., \acrshort{OCP} formats, \acrshort{IEEE} \acrshort{P3109}), vendor- or accelerator-specific designs (e.g., \acrshort{HiF8}, \acrshort{TSL} formats), and alternative representations (logarithmic formats, posits, and \acrshort{Takum}s). Table~\ref{tab:number_formats_desiderata} summarizes where these formats stand relative to the desiderata above, highlighting a key theme of the state of the art: \emph{hardware support is advancing faster than specification, portability, and analysis}.

Given the rapidly changing landscape of available formats,
which may be implemented differently by different vendors,
e.g. with different rounding modes, how exceptions are handled,
etc., we decided to fill in \cref{tab:number_formats_desiderata}
as follows. We base the properties described in each column 
either on
(1) a specification if all (or most) implementations adhere to that
specification, 
(2) a particular (widely available) implementation, or 
(3) a specification if an implementation is not yet available in hardware.
We document these choices below.

\subsubsection{Desiderata}
We begin by defining the desiderata that we will consider, 
i.e. the rows in \cref{tab:number_formats_desiderata}. Not all desiderata may be needed or desired for a particular application, so this should be treated as a menu of options to guide the choice of number format. 

\vspace{0.5em}
\noindent
{\bf Availability properties}

Unless otherwise specified, these properties in \cref{tab:number_formats_desiderata} are listed with yes (\yes), 
no (\no) or maybe (\maybe).

\begin{itemize}
    \item {\bf Hardware.} Is it commercially available (\yes), proprietary (\maybe), or not implemented in hardware (\no)? 
    \item {\bf Not Patented.} Is the use of the format encumbered by patents? This is a proxy for how widely the format can be expected to be supported.
    \item {\bf Portable.} Is it implemented consistently across platforms? \say{Maybe} (\maybe) will be used when the specification
    potentially allows many different implementations that are not all
    pairwise consistent. Still, since it is not implemented yet, like \acrshort{P3109},
    we cannot say for sure.
\end{itemize}
{\bf Design properties}
\begin{itemize}
    \item {\bf Well defined.} Are all aspects of the format well defined? This includes exception handling and propagation, as well as the accuracy of operations.
    \item {\bf Parameterized.} Is it defined for bit strings of variable lengths, or other parameters, such as the number of mantissa bits?
    \item {\bf Uniqueness.} Does every value (including exceptional values) have a unique bit representation? Or are many bit patterns potentially wasted? For example, \acrshort{IEEE} 754 includes
    many \acrshort{NaN}s, which are allowed but not required to
    be used to indicate uninitialized or missing 
    data (as in the R language) or include diagnostic
    information.
    \item {\bf Logarithmic format.} Is it designed 
    to represent numbers using their logarithms?
    We note that we do not include conventional formats with an
    exponent and 1 bit of precision, since only integer powers
    of 2 can be represented (e.g., Binary8p1se in \acrshort{P3109}).
    \item {\bf Saturation Modes Specified.} Can saturation behavior be specified independently of rounding modes? If saturation is not specified or if it is always on in a format, then it gets a \no.
    \item {\bf Rounding Modes Specified.}  Is the rounding of the results of operations specified? Does it include \acrfull{RTO}~\cite{bome08}) and \acrfull{SR}~\cite{cfhm22}? 
    (\acrlong{SR} means rounding up or down
    with higher probability to the closer floating
    point number.)
    Grade: 
    \begin{itemize}
    \item A: \acrfull{RN}, \acrfull{RU}, \acrfull{RD}, \acrfull{RZ}, \acrshort{SR}, and \acrshort{RTO} are specified.
    \item B: \acrshort{RN}, any of \acrshort{RU}, \acrshort{RD}, \acrshort{RZ}, and either \acrshort{SR} or \acrshort{RTO} are specified.
    \item C: \acrshort{RN}, \acrshort{RU}, \acrshort{RD}, and \acrshort{RZ} are specified.
    \item D: \acrshort{RN} and either \acrshort{SR} or \acrshort{RTO} are specified.
    \item E: Only \acrshort{RN} is specified.
    \item F: No rounding specification is provided.
    \end{itemize}
    \item {\bf Block formats.} Does the format support a 
    block format, i.e. an
    array of numbers in a common format along with
    at least one scaling factor shared by all the numbers,
    yes (\maybe) or no (\no)? Is it well-defined in terms of edge cases (\yes)?
    \item {\bf Format Convertible.} Is conversion between formats of the same kind (i.e., same column in 
    \cref{tab:number_formats_desiderata}) with differing parameters straightforward? For example, converting
    from a format with a \acrshort{NaN} to one without does not always work.
    \item {\bf Precision Convertible.} Is conversion between bit strings of different lengths, i.e. precisions, specified? 
    For example, it might
    be difficult to round from one logarithmic
    format to a shorter one correctly.
    Not applicable (\maybe) for formats that have only one precision setting. 
    \item {\bf High precision accumulator.} Does this format include a high(er) precision accumulator in at least one implementation of dot product or matrix multiply from a major vendor: yes (\maybe) or no (\no)? And if yes, is it accessible (\yes) to the user?
\end{itemize}

{\bf Numerical properties}
\begin{itemize}
    \item {\bf Numerical analysis exists or is possible.} Does existing numerical analysis apply immediately, at least for
    small enough problems (\yes)? If not, are there works on numerical analysis in the format, yes (\maybe) or no (\no)?~\cite{demm87}
    \item {\bf Composable for emulation of higher precision.}
    Can a format be used to emulate higher precision numbers, and (some) operations on them, as described in \cref{sec:emulation}? \yes{} if algorithms and experiments have been shown in literature or implementations, and \maybe{} if no algorithms for emulation in the format exist at the time of writing. Note, algorithms for emulation with \acrshort{IEEE} 754-like formats have been extensively explored, so all formats that encode numbers like \acrshort{IEEE} 754 get a \yes{}.

    \item $\mathbf{a-b=0 \Leftrightarrow  a=b}$. Does this elementary property hold for finite $\mathbf{a}$ and $\mathbf{b}$? This typically requires support for subnormal numbers (i.e. numbers with the smallest
    exponent that can have leading zeros in their 
    mantissas)
    without which it breaks if $\mathbf{a}$ and $\mathbf{b}$ are the two smallest consecutive non-zero representable numbers.
    \item $\mathbf{e_{\textrm{max}}+e_{\textrm{min}}}$: Take $e_{max}$ and $e_{min}$ to be maximum and minimum exponents in a format, respectively, and take $e_{\textrm{max}}+e_{\textrm{min}}$. Is this sum close to zero, i.e., is there symmetry in the dynamic range of the representation?
\end{itemize}

\subsubsection{Formats} We consider the desiderata for the following list of formats, which are the most widely discussed but certainly not exhaustive.

\vspace{0.5em}
\noindent
{\bf \acrshort{IEEE} 754 and derivatives}
\begin{itemize}
    \item{\bf 754} \acrshort{IEEE} 754 standard~\cite{ieee19}.
    \item{\bf \acrshort{BF16}} \acrlong{BF16}. Derived from \acrshort{FP32} by cutting the mantissa from 24 to 8 bits (including the implicit leading 1)~\cite{inte18,BFloat16}. 
    \item{\bf \acrshort{TF32}} NVIDIA's \acrlong{TF32} (19 bits arithmetic, 32 bits storage)~\cite{stmi21}. Equivalent to \acrshort{FP16} mantissa plus \acrshort{FP32} exponent.
\end{itemize}

{\bf \acrshort{AI}-oriented proposals}
\begin{itemize}
    \item{\bf \acrshort{OCP}} \acrlong{OCP} \enquote{standard} for OFP8 (scalar)~\cite{ocp23a} and OFP-MX (block)~\cite{ocp23b}.
    \item{\bf \acrshort{P3109}} A proposed standard, \acrshort{IEEE} \acrshort{P3109}, for floating-point formats and arithmetic for machine learning~\cite{ieee25}.
\end{itemize}
{\bf Vendor formats}
\begin{itemize}
    \item{\bf \acrshort{AGQ}} Formats without infinities, minus zero, and with only one \acrshort{NaN}, proposed by \acrshort{AMD}-GraphCore-QualComm~\cite{njjm22}\cite{AGQ}.
    May have been superseded by \acrshort{OCP} formats, since \acrshort{AMD} uses \acrshort{AGQ} in the CDNA3 architecture and \acrshort{OCP} in the CDNA4 architecture \cite{CDNA4}. 
    \item{\bf \acrshort{HiF8}} HiFloat8 format proposed by Huawei~\cite{lzwl24}. 
    \item{\bf \acrshort{TSL}} Tesla Dojo specific format~\cite{tesla}.
    \item{\bf \acrshort{LNSM}} NVIDIA's logarithmic number format used in their LNS-Madam paper~\cite{zdvz22}. 
    Tensordyne \cite{Tensordyne} also uses a logarithmic number format for its \acrshort{AI} chip, but aside from being 16-bit, little has been disclosed about this format. 
\end{itemize}
{\bf Alternative representations}
\begin{itemize}
    \item{\bf Posit} A format with ``tapered precision'' that has higher precision for numbers with
    magnitude near 1, falling off as the exponent gets 
    further from 0, i.e., as the magnitude grows or shrinks.  \cite{guyo17, posi22} 
    \item{\bf \acrshort{Takum}} Constant dynamic range format tailored towards scientific computing within the posit format family~\cite{hunh24}.
\end{itemize}

\begin{table*}[tbp]
\centering
\caption{Desiderata for floating-point formats.}
\label{tab:number_formats_desiderata}
\footnotesize
\bgroup
\renewcommand{\arraystretch}{1.25}
\begin{tabular}{|p{3.1cm}|c|c|c|c|c|c|c|c|c|c|c|}
  \hline
                                     & \textbf{754} & \textbf{BF16} & \textbf{TF32} & \textbf{OCP} & \textbf{P3109} & \textbf{AGQ} & \textbf{HiF8} & \textbf{TSL} & \textbf{LNSM} & \textbf{Posit} & \textbf{Takum} \\ \hline
Hardware                             & \yes         & \yes          & \yes          & \yes         & \no            & \yes         & \maybe        & \yes         & \no        & \yes           & \no            \\ \hline
Not Patented                         & \yes         & \yes          & \yes          & \yes         & \yes           & \yes         & \no           & \yes         & \no        & \yes           & \yes           \\ \hline
Portable                             & \yes         & \no           & \no           & \no          & \maybe         & \yes         & \maybe        & \no          & \maybe     & \no           & \maybe           \\ \hline\hline
Well Defined                         & \yes         & \no           & \no           & \no          & \yes           & \no          & \no           & \no          & \no        & \yes           & \yes           \\ \hline
Parameterized                        & \yes         & \no           & \no           & \no          & \yes           & \no          & \yes          & \no          & \yes       & \yes           & \yes           \\ \hline
Uniqueness                           & \no          & \no           & \no           & \no          & \yes           & \yes         & \yes          & \yes         & \yes       & \yes           & \yes           \\ \hline
Logarithmic format                   & \no          & \no           & \no           & \no          & \no            & \no          & \no           & \no          & \yes       & \no            & \yes           \\ \hline
Saturation specified           & \no          & \no           & \no           & \yes         & \yes           & \yes         & \yes          & \no          & \no        & \no           & \no           \\ \hline
Rounding specified             & C            & F             & F             & E            & A              & E            & D             & D            & E          & E         & E         \\ \hline
Block formats                        & \no          & \no           & \no           & \maybe       & \yes           & \no          & \no           & \no          & \no        & \no            & \no            \\ \hline
Format Convertible                   & \yes         & \yes          & \yes          & \yes         & \yes           & \yes         & \yes          & \yes         & \yes       & \yes           & \yes           \\ \hline
Precision Convertible                & \yes         & \maybe        & \maybe        & \no          & \yes           & \no         & \no        & \no          & \no       & \yes           & \no           \\ \hline
High prec. accumulator           & \yes         & \yes          & \yes          & \yes         & \maybe         & \yes         & \no           & \no          & \yes       & \yes           & \no           \\ \hline\hline
Numerical Analysis                   & \yes         & \yes          & \yes          & \yes         & \yes           & \yes         & \no           & \yes         & \no        & \maybe         & \maybe         \\ \hline
Composable for emul.             & \yes         & \yes          & \yes          & \yes         & \yes           & \yes         & \yes          & \yes         & \maybe        & \maybe            & \maybe         \\ \hline
$a=b \Leftrightarrow a-b=0$          & \yes         & \maybe        & \maybe        & \yes         & \yes           & \yes         & \no           & \yes         & \no        & \yes           & \yes           \\ \hline
$e_{\textrm{max}}+e_{\textrm{min}}$  & 1            & 1             & 1             & \{1, 2\}     & 0              & 0            & 0             & Varied       & Varied     & 0              & 0              \\ \hline
\end{tabular}
\egroup
\end{table*}

\subsection{Future directions}
It can be expected that the Cambrian explosion of number formats we witnessed in the last decade will settle down into a funneling toward fewer, more mature number formats. This is further accelerated by the fact that hardware will remain relatively short-lived, and advances in shorter tape-out turnaround times will remain part of the chip industry, allowing it to better react to such trends in number formats. Unlike \acrshort{AI}, which mostly relies on \acrshort{FMA} and specific hardware-based analytical functions, \acrshort{HPC} relies more on high-accuracy operations and functions, and the relative dominance of \acrshort{AI} will dictate the availability of such operations in hardware.
\par
While hierarchical data formats (\enquote{block scaling}) see widespread use in \acrshort{AI},
they are still a niché application in \acrshort{HPC}. It remains to be seen how the philosophical
shift away from single-number machine number representations to multi-number hierarchical
representations is reflected in number format design and adaptable to \acrshort{HPC} workloads,
as it requires regularity in the underlying datasets. From the current state of tentative block scaling in \acrshort{P3109}, we will see more such influences on number format design. The term \enquote{number format} will, in the same vein, reflect not only the mapping of a given number of bits to a represented value, but also the hierarchy of the surrounding data.
\par
We will hopefully see a move away from simply stating flop/s, given there is a significant difference between a 4- or 64-bit floating-point
operation. 
What truly matters is the amount of information that is processed in a given amount of time and amount of energy.
This can be simply reflected at the bit level using entropic evaluation \cite{havu2025}, which allows distinguishing between different bit counts. However, with floating-point formats, one also has to account for differences among formats in such a comparison.

\subsection{Open questions}
One main set of open questions concerns emulation (see \cref{sec:emulation}). While addition and multiplication are handled in the current emulation schemes properly, it is necessary for many applications to support other operations, like
division, square root, and various transcendental functions. 
A question remains as to how
well this can be done with such emulation schemes, including the question
of how one can do proper exception handling. The overall question emerges if
it is better to have hardware/software support for many low- and high-precision formats,
or if one builds everything, using emulation, on top of one to two low-precision formats. The main aspect of differentiation is that operations in lower precisions
are orders of magnitude faster and cache-local; however, compilers on a higher level
could be better than a hardware implementation, as they can pipeline these operations
better due to a better overview of the computational flow.

\par
Another question is the relative importance of operations. Depending on applications,
some operations could have a higher weight compared to others. With the increasing use of tailored
hardware, what consequences arise from manufacturers deliberately omitting standardized
operations, for instance, division or square root? To further assess this aspect, more
empirical data needs to be collected on how different operations are used in respective
applications. The same applies to number distributions, which dictate the choice of
formats for specific applications. While this is more and more studied for \acrshort{AI}, it is
an open research question in \acrshort{HPC}.

\ignore{
\subsection{Future directions}
It can be expected that the Cambrian explosion of number formats we witnessed in the last decade will settle down into a funneling toward fewer, more mature number formats. This is further accelerated by the fact that hardware will remain relatively short-lived, and advances in shorter tape-out turnaround times will remain part of the chip industry, allowing it to better react to such trends in number formats. Unlike \acrshort{AI}, which mostly relies on FMA and specific hardware-based analytical functions, \acrshort{HPC} relies more on high-accuracy operations and functions, and the relative dominance of \acrshort{AI} will dictate the availability of such operations in hardware.
\par
While hierarchical data formats (\enquote{block scaling}) see widespread use in \acrshort{AI},
they are still a niché application in \acrshort{HPC}. It remains to be seen how the philosophical
shift away from single-number machine number representations to multi-number hierarchical
representations is reflected in number format design and adaptable to \acrshort{HPC} workloads,
as it requires regularity in the underlying datasets. From the current state of tentative block scaling in \acrshort{P3109}, we will see more such influences on number format design. The term \enquote{number format} will, in the same vein, reflect not only the mapping of a given number of bits to a represented value, but also the hierarchy of the surrounding data.
\par
We will hopefully see a move away from simply stating FLOPS (floating-point operations
per second), given there is a significant difference between a 4- or 64-bit floating-point
operation. What truly matters is the amount of information that is processed in a given amount of time. This can be simply reflected at the bit level using entropic evaluation \cite{havu2025}, which allows distinguishing between different bit counts. However, with floating-point formats, one also has to account for differences among formats in such a comparison.
} %

\section{Floating-point emulation}\label{sec:emulation}

In recent years, the development of \acrshort{AI} algorithms has increased \acrshort{GPU} usage and led to the emergence of various low-precision formats with 16, 8, or fewer bits. The performance gap between \acrshort{FP64} arithmetic and these lower-precision is becoming increasingly favorable to lower-precision formats, particularly when using matrix compute units. Many scientific applications, however, require the accuracy provided by \acrshort{FP64} or higher-precision formats.

When not available in hardware, high precision can be emulated using software approaches that can deliver both accuracy and performance. Floating-point emulation consists of representing a high-precision number as a structured combination of several lower-precision values (floating-point or integer), and performing operations on this representation so as to approximate higher-precision arithmetic while exploiting the speed of low-precision hardware.  This approach can deliver both accuracy and performance: emulation increases the number of arithmetic operations, but this cost can be offset by performing lower-precision computations.

Emulation provides access to high precision on hardware where native support is lacking, and allows scientific applications to fully exploit fast, low-precision units while still meeting accuracy requirements. It can therefore enable substantial gains in time and energy to solution, especially in compute-bound settings.

It is important to note, however, that emulated arithmetic is not bit-wise equivalent to native high precision. In particular, the dynamic range can be limited and the rounding behavior may differ.

\subsection{State of the art}

Floating-point emulation techniques can be broadly divided into two families: one based on multiword arithmetic and one based on the Ozaki scheme. In the literature, both approaches have primarily been studied in the context of compute-bound kernels, with a particular emphasis on matrix multiplication.

\subsubsection{Multiword arithmetic} Multiword arithmetic is a well-established technique to enable extended-precision computations \cite[Chap. 14]{mbdj18}. Initially designed for scalar operations, this technique can also be used to emulate high-precision matrix multiplication with lower-precision operands. Traditionally, multiword algorithms rely on error-free transformations, but in the case of matrix multiplication, a similar result can be obtained by using high-precision accumulation of dot products. This technique can only be used efficiently to simulate up to the accumulator's precision.

\begin{figure*}
    \begin{center}
\includegraphics[width=.7\linewidth]{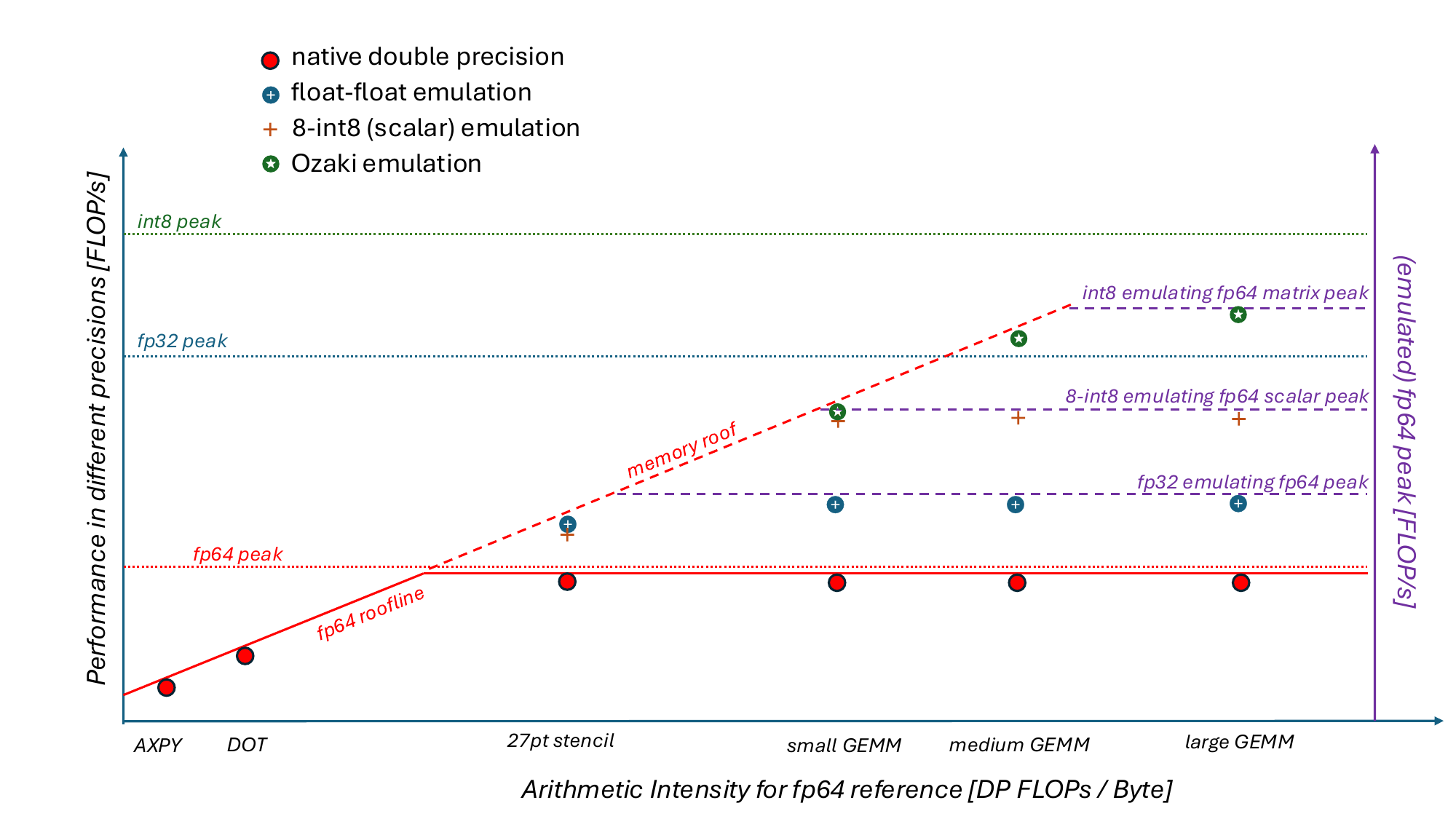}
    \caption{Performance of different kernels using native FP64 or emulation.}
    \label{fig:roofline}
    \end{center}
\end{figure*}

This emulation strategy does not yield the same result as the native high-precision arithmetic. First, the range of the emulated format is restricted to the range of the low-precision format used in the emulation---for example, when emulating double precision using single-precision arithmetic, the range of the emulated format is the range of single precision. Second, the combined significands of two single-precision values (24 significand bits each) fall 5 bits short when emulating double precision.

The traditional multiword emulation on a scalar level is flexible and can be used for both memory-bound and compute-bound algorithms. However, a compute-bound algorithm can become memory-bound if the emulation is based on a format that has significantly higher peak performance on the underlying hardware architecture, as shown in \cref{fig:roofline}. For compute-intensive algorithms, the scalar arithmetic units can become the performance-limiting factor when using low-precision formats in the emulation.

The emergence of high-throughput mixed-precision matrix multiplication units, such as NVIDIA Tensor Cores and \acrshort{AMD} Matrix Cores that take \acrshort{FP16}, \acrshort{BF16}, or \acrshort{TF32} inputs and perform accumulation in \acrshort{FP32}, has opened new opportunities to exploit their computational performance beyond their original purpose.
Several authors have used these units to emulate matrix multiplication at the precision of the accumulator \cite{markidis_nvidia_2018,fhlm23,ootomo_recovering_2022}.
Such functionality is already partially available in libraries such as cuBLAS\footnote{\url{https://developer.nvidia.com/cublas/}} and CUTLASS,\footnote{\url{https://docs.nvidia.com/cutlass}} which support mixed-precision \acrshort{GEMM} with \acrshort{FP32} accumulation.

\subsubsection{Ozaki-style algorithms}  More recently, the integer matrix multiplication units available on \acrshort{GPU}s have renewed interest in the emulation strategy known as the ``Ozaki'' scheme.
The high-level idea is to emulate \acrshort{FP32} or \acrshort{FP64} matrix multiplication using integer Tensor Cores, which accumulate the product of 8-bit matrices (\acrshort{INT8}) in 32-bit integer arithmetic (\acrshort{INT32}) \cite{ootomo_dgemm_2023}. In some detail, the input matrices $\mathbf{A}$ and $\mathbf{B}$ are scaled, decomposed into matrix slices, and processed as follows.
\begin{enumerate}
  \item Normalization and Error-Free Splitting: First, each row of $\mathbf{A}$ and each column of $\mathbf{B}$ are normalized (scaled) by their respective absolute maximum values. After scaling, the matrices are decomposed into slices $\mathbf{A}=\sum_i \mathbf{A}^{(i)}$ and $\mathbf{B}=\sum_j \mathbf{B}^{(j)}$ so that the product $\mathbf{A}^{(i)} \cdot \mathbf{B}^{(j)}$ does not overflow the \acrshort{INT32} accumulator.
  \item Each product $\mathbf{A}^{(i)} \cdot \mathbf{B}^{(j)}$ is evaluated using \acrshort{INT8} Tensor Cores.
  \item Accumulate the resulting matrices in double precision: The resulting \acrshort{INT32} products are scaled back to their original floating-point magnitudes by considering both the splitting factors and the initial normalization factors. These values are then accumulated into an \acrshort{FP64} buffer. 
\end{enumerate}
By combining these exact partial products, the final result typically achieves better-than-\acrshort{FP64} accuracy, effectively emulating high-precision arithmetic on low-precision hardware.
Besides \acrshort{INT8} Tensor Cores, the Ozaki scheme can also operate on \acrshort{FP16} \cite{mukunoki_reproducible_2020} or FP8 Tensor Cores \cite{mukunoki_dgemm_2026}. Note that in error-free splitting, the number of significand bits in each slice depends on the accumulator length. On NVIDIA \acrshort{GPU}s, double-precision \acrshort{GEMM} using the Ozaki scheme is readily available in NVIDIA cuBLAS, making it accessible to anyone. Furthermore, dense linear algebra algorithms that rely on the Ozaki scheme for \acrshort{GEMM} as the central building block have been developed in double precision.
For sufficiently large matrices, the performance of dense linear algebra operations based on the Ozaki scheme is limited by the performance of the \acrshort{INT8} Tensor Cores. For smaller matrices, e.g., handled in a batched fashion, the matrix multiplication using the Ozaki scheme is expected to become memory-bound, as shown in~\Cref {fig:roofline}.
In particular, in those scenarios, the current workflow, which handles matrix splitting in main memory and thereby increases communication, is not ideal. Implementations targeting small sizes should ideally handle the splitting and all slice access in fast on-chip memory.
The Ozaki Scheme II is an emulation method that performs significand decomposition based on the Chinese Remainder Theorem \cite{ozaki2025ozakischemeiigemmoriented,uchino_oz2_2025}.
Since the computational cost scales linearly with the number of slices, the method is particularly efficient for large matrix sizes, where the decomposition overhead is effectively amortized.

Table~\ref{tab:Emulation_perf} presents for 
matrix multiplication carried out using 
multiword arithmetic or fixed point emulation using Ozaki scheme:  
the number of matrix multiplications,
the number of floating-point operations, 
the performance with respect to \acrshort{FP64} of  \acrshort{FP16} on H100 NVIDIA \acrshort{GPU} and \acrshort{INT8} on B200 NVIDIA \acrshort{GPU}. For a given hardware, this allows us to assess whether the Ozaki scheme offers performance benefits over native double precision.

\begin{table*}
\centering
\caption{Emulation performance for matrix multiplication 
 ($n$ denotes the matrix size, 
 $p$ the number of floating-point numbers in the multiword structure,
 $s$ the number of slices in the Ozaki scheme.)
 }
\begin{tabular}{lcccc}
\toprule
& \#MM  & flops & FP64/FP16 (H100) & FP64/INT8 (B200)\\
\midrule
 High precision  & 1 &  $\approx 2n^3$ & 1 & 1\\
 Multiword MMA based        & $p(p+1)/2$ &  $\approx p^2n^3$ & 5  &  20\\ 
 Fixed-point emulation, Ozaki I& $s(s+1)/2$ &  $\approx s^2n^3$ &8 & 32 \\ 
 Fixed-point emulation, Ozaki II& $s$ &  $\approx 2s n^3$  & 30 & 120\\
 \bottomrule
 \end{tabular}
 \label{tab:Emulation_perf}
\end{table*}

\subsection{Future directions}

The techniques discussed in this section can emulate any precision, and while it is likely that the main focus of research will remain on \acrshort{FP64}, there will be a need to develop emulation techniques for \acrshort{FP32} and for higher-than-64-bit formats.

Drop-in replacements for the most popular dense \acrshort{BLAS} Level 3 operations exist under the Ozaki scheme. For Level 1 and Level 2 operations and sparse linear algebra, such drop-in replacements would not be attractive in most cases, because the kernels are memory-bound, and kernel fusion strategies that reduce memory access are more useful in this case.
Therefore, for operations with low arithmetic intensity, it would be more productive to develop codelets based on emulation that can be combined, for example, by the compiler, to form complex kernels.

To inform hardware manufacturers' decisions, it may be useful to provide a pen-and-paper calculation to assess which features hardware should provide for emulation to be beneficial. In particular, one should understand what combinations of precisions should be available, and how much faster than mono-precision units these mixed-precision units should be. Preliminary calculations in this direction are shown in \cref{tab:Emulation_perf}, which presents the minimum performance gap required to make emulation faster than native \acrshort{FP64} precision on NVIDIA H100 and B200 \acrshort{GPU}s, based on the vendor's declared theoretical peak performance.

%

\subsection{Open questions}

The numerical semantics of emulation remain an open and central set of questions.
A key issue is identifying the optimal balance between the exponent range and the precision of the significand.
It is not clear whether faithfully reproducing the full dynamic range of \acrshort{FP64} is necessary, or whether maintaining higher significand precision would be more useful in practice.
Closely related is the treatment of special floating-point values, such as infinities, negative zero, and \acrshort{NaN}s.
It remains unclear how these can be detected, propagated, and interpreted in emulation, especially in mixed-precision contexts.
Beyond standard floating-point representations, an important direction is the extension of emulation techniques to alternative numerical formats and guarantees. This includes the efficient emulation of emerging formats such as Posit and \acrshort{Takum} using low-precision hardware, as well as the application of emulation to interval arithmetic, where the goal is to combine high performance with mathematically certified bounds for critical computations.

A second group of questions concerns portability and reproducibility.
Whether emulation should aim for cross-vendor reproducibility, and if so, how such reproducibility should be defined, remains unclear.
Differences in \acrshort{GPU} architectures, low-precision units, and optimization strategies across vendors and hardware generations complicate the notion of consistent behavior. This raises fundamental questions about the trade-offs between performance, hardware-specific optimization, and numerical consistency, as well as whether a universally reproducible emulation layer is feasible---or even desirable---in practice.

Finally, there are significant challenges in automation and system support for emulation.
One major goal is to design a compiler or runtime system capable of automatically identifying where emulation is beneficial, inserting it selectively and managing it efficiently without burdening the programmer.
It is still an open question whether this should be addressed at the library, compiler, or language level.
Hardware may be able to actively assist software emulation, providing hints or error metrics to make automatic, workload-driven emulation practical and efficient.
Future hardware could expose fine-grained error metrics, dynamic-range monitors, or precision-confidence signals that help the runtime to promote operations to higher precision when necessary.

Emulation is one path to trusted low precision, but whether it is worthwhile depends on the hardware trends (\cref{sec:emerging}), the semantics of the number formats available (\cref{sec:numberformats}), and on the software/runtime support (\cref{sec:softwaredesign} and \cref{sec:multilevel}).

\section{Emerging architectures}
\label{sec:emerging}

Emerging computer architectures are changing the role of numerical precision in scientific computing. 
For decades, reliability in computing has been supported through hardware
mechanisms that detect and correct physical faults (e.g., ECC for memory).
In scientific computing, \acrshort{FP64} has played an
analogous role for \emph{numerical reliability}. 
Approximation and rounding errors are inevitable
in finite-precision arithmetic, 
but \acrshort{FP64} is often sufficient so that software developers do not need to be explicitly concerned about them~\cite{higham2002accuracy}.
That margin is now shrinking as modern systems increasingly optimize for energy efficiency, throughput, and data movement rather than universally high precision.

Three architectural trends drive this transition. First, the power wall and the growing dominance of data movement mean that moving data through the memory hierarchy often costs far more—in both energy and time—than the arithmetic itself. As a result, reducing operand size and communication volume has become a first-order optimization target. Second, hardware is rapidly shifting toward lower-precision execution units, driven largely by \acrshort{AI} workloads, while native \acrshort{FP64} increasingly becomes a secondary capability on many accelerators. Third, architectures are beginning to incorporate reliability-oriented numerical features such as \acrshort{SR}, wide accumulators, fused and deferred rounding schemes, and mixed-precision matrix engines. These mechanisms suggest a future in which precision is no longer a fixed global setting, but a dynamically managed resource.

The next stage of this evolution is not simply ``more low precision.” Instead, architectures may provide precision-aware hardware services that actively participate in numerical reliability. Examples include local precision-risk indicators within \acrfull{FPUs}, hardware support for hybrid correctness loops that combine approximate execution with selective certification, and cross-architecture abstractions that expose approximate and corrective operators independently of device-specific implementations. Such mechanisms could allow systems to exploit aggressive approximation while still maintaining scientifically meaningful correctness guarantees.

This shift raises several open questions. At the hardware level, what lightweight signals can expose local numerical risk without excessive overhead? At the system level, how should runtimes and compilers use these signals to trigger adaptive precision refinement, recomputation, or fail-fast behavior? At the validation level, how should \acrshort{VaV}, and \acrshort{UQ} evolve when hardware-induced numerical variability becomes part of the computational model itself?

Ultimately, emerging architectures can expose new precision opportunities, but realizing them safely and portably depends on the analyzable number formats of \cref{sec:NumberFormats} and the hardware/software abstractions developed in \cref{sec:co-design}.

\paragraph*{From numerical ``safety margins'' to precision-aware hardware.}
As the ecosystem moves toward lower-precision
formats, the risk of roundoff error, cancellation, and numerical instability
increases, and we need mechanisms that can indicate when reduced precision is
safe --- and when it is not. 
At the same time, the dominant energy and performance
cost in modern systems is increasingly \emph{data movement}, not arithmetic.
Moving wide operands through the memory hierarchy can dwarf the cost of the
associated floating-point operations~\cite{dally2022model,DataMovement}.
This asymmetry suggests an architectural opportunity. Modest additional logic
in \acrfull{ALUs}/\acrshort{FPU}s may be comparatively low-cost, yet valuable, if they can provide
runtime indicators of precision loss or rounding sensitivity.

One concrete direction is \emph{hardware-assisted numerical introspection}.
For example, recomputing an operation under two rounding modes and comparing
results can provide a local estimate of rounding sensitivity, serving as an
\say{early warning} signal for precision loss. Conceptually, such indicators
play a role akin to parity in fault detection: they do not solve correctness
end-to-end, but they can flag cases that merit correction, refinement, or
fallback. The broader research question is how to evolve \acrshort{ALU}s/\acrshort{FPU}s from passive arithmetic engines into active contributors to numerical reliability---through
lightweight support for emulation, rounding sensitivity tracking, and low-cost
error proxies---without compromising performance.

\subsection{State of the art}

\subsubsection{Energy and the power wall.}
\begin{figure*}
    \centering
    \includegraphics[width=.7\linewidth]{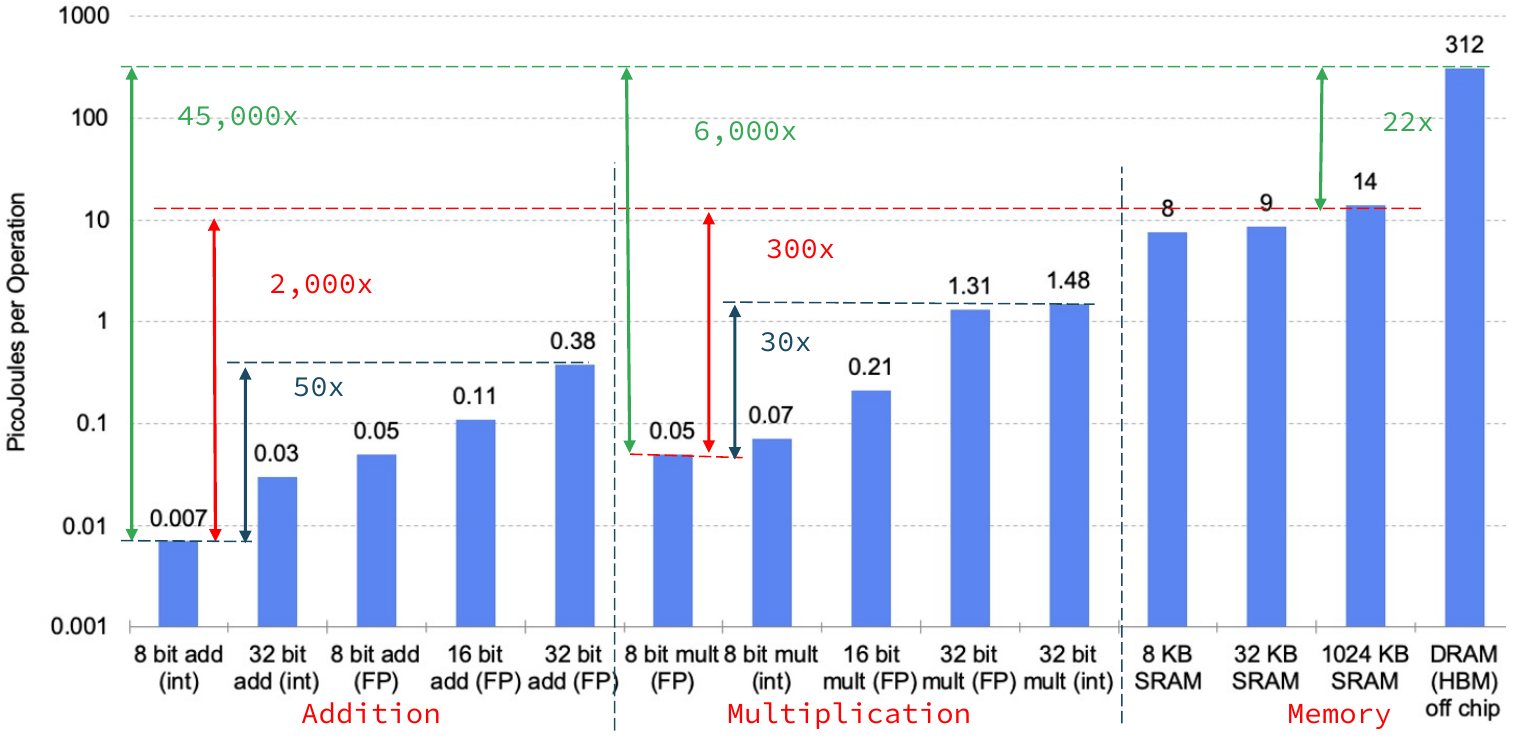}
    \caption{Energy costs of memory operations and floating-point operations.
    Figure by David Patterson (Google). 
    }
    \label{fig:Patterson}
\end{figure*}

The widening gap between arithmetic throughput and memory bandwidth is reflected
not only in performance but also in energy: data movement dominates power
consumption in contemporary architectures~\cite{dally2022model,DataMovement}.
Some estimates (\cref{fig:Patterson}) suggest orders-of-magnitude energy
differences between simple integer operations and off-chip memory accesses.
Energy breakdowns likewise attribute only a tiny fraction of total energy to the
\acrshort{FPU} itself, with most energy spent moving data through caches, registers,
and interconnects (e.g. the datapath) as shown clearly in  \cref{fig:Patterson} from Patterson.

\subsubsection{Industry/Market shift toward lower precision}
Driven largely by \acrshort{AI}, \acrshort{FP64} and \acrshort{FP32} are increasingly complemented by compact
formats (8, 6, or even 4 bits) as described in \cref{sec:emulation}.
\acrshort{GPU}s now provide extensive support for \acrshort{FP8} and
below, often alongside mixed-precision execution. \acrshort{AMD} similarly integrates
arithmetic for low-precision floating-point and integer types, including
micro-scaled formats aligned with \acrshort{MXP} variants~\cite{ocp_mx_spec_2023}.
\acrshort{CPU}s are following more slowly: Arm increasingly supports \acrshort{FP16}/\acrshort{BF16} and \acrshort{INT8}
acceleration via \acrshort{SIMD}/vector extensions, while RISC-V is evolving optional
extensions for \acrshort{BF16}/\acrshort{FP16} and is sometimes paired with low-precision accelerators. 

\subsubsection{Relative marginalization of native FP64}
Also, as described in \cref{sec:emulation}, as \acrshort{GPU}s 
devote more area to low-precision matrix engines (e.g., Tensor Cores),
native \acrshort{FP64} capability risks becoming a secondary design point for mainstream
accelerators. One response is \emph{emulation}: for high–arithmetic-intensity
kernels such as \acrshort{DGEMM}, schemes such as Ozaki-style approaches can emulate \acrshort{FP64}
(or other precision) using integer arithmetic, sometimes delivering high
effective throughput on hardware optimized for low-precision. However, such
approaches can require sophisticated implementations and incur nontrivial
computational and memory overheads; they are most compelling for large dense
kernels and are far less attractive for memory-bound computations or on \acrshort{CPU}s
without massive low-precision matrix engines.

\subsubsection{Emerging number formats}
\label{subsec:EmergingNumberFormats}
Beyond \acrshort{IEEE}~754, alternative formats aim to improve accuracy-per-bit or simplify
certain operations as described in more detail in \cref{sec:NumberFormats}. 
For example, Posits provide tapered precision and can use wide accumulators
to improve dot-product accuracy~\cite{guyo17,gustafson2019end}.
\acrshort{Takum}s pursue different tradeoffs and have been proposed as simpler to implement
in hardware for scientific computing~\cite{hunh24,hunhold2024codec,hunhold2025simd}.
They offer a broad dynamic range at low bit-widths and can change mixed-precision
tradeoffs by decoupling some accuracy considerations from dynamic-range concerns.
At present, both Posits and \acrshort{Takum}s remain primarily research efforts, implemented
mostly via \acrshort{FPGA}s/RISC-V prototypes or experimental \acrshort{ASIC}s, with limited commercial
adoption.

\subsubsection{Stochastic rounding}
\acrshort{SR} rounds probabilistically to a neighboring representable
values, reducing systematic bias that can accumulate under deterministic rounding,
especially at low-precision~\cite{chm21,cfhm22}. Despite its algorithmic value,
\acrshort{SR} is not yet a mainstream \acrshort{ALU} feature and is primarily found in selected,
\acrshort{AI}-focused accelerators.

\subsubsection{Rounding-deferred summation and wide accumulators}
\acrshort{RDS} accumulates intermediate results in higher or
exact precision and rounds only once at the end, reducing error accumulation in
dot products and related kernels~\cite{higham2002accuracy}. Classic concepts
include Kulisch-style accumulators~\cite{kulisch2002computer} and the Posit 
\emph{quire}. While \acrshort{FP64}-exact accumulation may require very wide registers,
\acrshort{RDS} becomes much more plausible for sub-16-bit formats, where a compact
accumulator can capture the full dynamic range. Current hardware provides only
partial support via \acrshort{FMA}, which reduces (but does not
eliminate) intermediate rounding. More complete support for wide or exact
accumulation would remove a major source of numerical error in contemporary
reductions and dot products and aligns naturally with mixed-precision \acrshort{MMA} units~\cite{bhlm20}.

\subsection{Future directions}
\label{future}

\subsubsection{Algorithm--architecture co-design for precision}
When an application is memory-bound, additional \acrshort{FP64} units do not increase
performance; data movement dominates. In that regime, native \acrshort{FP64} may be less
efficient than emulation strategies that exploit high-throughput low-precision
paths. For compute-bound kernels (e.g., \acrshort{DGEMM}), emulating \acrshort{FP64} atop tensor
hardware can be attractive, but making emulation broadly useful will require
architectural support to reduce fixed overheads (e.g., \acrshort{ISA} support for double primitives, extended accumulators, or other microarchitectural
assist). If emulation becomes fast and certifiable, the need for large native
\acrshort{FP64} resources may diminish for some workloads. These benefits, however, do not
carry over automatically to memory-bound \gls{BLAS}-1/2 style operations.

\subsubsection{Precision awareness at the ALU/FPU level}
As some correctness goals shift from strict bitwise reproducibility toward
statistical validation of quantities of interest, we need ways to identify
numerically risky regions without relying solely on global error estimators.
A promising direction is \acrshort{FPU}-local precision monitoring. Indicators of severe cancellation or rounding sensitivity could serve as a low-cost warning during
low-precision execution. Such signals are inherently local and cannot replace
algorithmic error analysis, but they could trigger selective recomputation,
adaptive precision refinement, or fail-fast behavior before instability spreads.

\subsubsection{Noise-tolerant phases and structured approximation}
Many workflows contain phases where full precision is unnecessary (coarse models,
preconditioners, Monte Carlo proposal steps, or \acrshort{AI} surrogates). Here, gains may
come less from \say{adding bits} and more from exploiting structure. Reloadable
lookup tables~\cite{ramanathan2020lutenergy,sutradhar2023flutpim,lin2024lutin}
and learned surrogates~\cite{diaw2024surrogates} can replace expensive floating
point with table- or model-based evaluation that respects known invariants.
This reframes the objective: match the representation to the \emph{effective
complexity} of the phase, allocating precision and computation where the problem
demands it rather than enforcing worst-case guarantees uniformly.  The potential for
this approach is described in more detail in \cref{sec:co-design}.

\subsubsection{Hybrid correctness loops: generation and certification}
Correctness can be framed as a two-stage loop. Use low-precision to generate
candidate solutions quickly, then certify (or refine) using higher precision,
residual checks, sensitivity analysis, or bounded error tests. As arithmetic
becomes cheaper relative to data movement, recomputation and multi-candidate generation may become economically viable when data remain local in the pipeline.
Realizing this pattern robustly will require hardware and compiler/runtime
support for rapid candidate evaluation, selective refinement, and efficient
fallback.

\subsubsection{Cross-architecture abstractions}
Cross-architecture abstractions can treat heterogeneous systems as compositions
of operators rather than device-specific implementations: (i) approximate
operators optimized for performance/energy with controlled inaccuracy, and
(ii) correction/certification operators that enforce numerical guarantees
or bounds~\cite{ApproxComputSurveyII}. If such interfaces can expose
architecture-independent \emph{guarantees} (analogous to how \acrshort{BLAS} exposes
architecture-independent \emph{semantics}), then precision management, emulation,
and device-specific optimizations can remain encapsulated beneath certified
operator contracts, enabling portability without sacrificing trust.

\subsection{Open questions}
\label{sec:mrgopen}

\subsubsection{Local dynamic precision signals}
What lightweight features can \acrshort{FPU}s provide to cheaply estimate local precision
loss and accumulated rounding sensitivity, while remaining method-agnostic?
Dual-rounding-mode evaluation is one example. The approach is to compute two
results, detect divergence, and then raise a flag/exception that triggers corrective action. In contrast, full interval arithmetic doubles datapath width and is
likely infeasible at scale; is there a reduced, \acrshort{FPU}-local variant that captures
some of its benefits without prohibitive data movement?

\subsubsection{Connecting local indicators to global error control}
\acrshort{FPU}-local indicators may be too myopic to guide algorithm-level correctness on
their own. Can we relate local signals to global error bounds and estimators?
For example, adaptive precision strategies may rely on norms, conditioning
estimates, or other global diagnostics to identify
numerically sensitive components. Can this be generalized across algorithm
families, including both residual-based methods (e.g., linear algebra) and
residual-free settings (e.g., \acrshort{UQ}, data assimilation)?

\subsubsection{Speculative execution with numerical commit/rollback}
Modern \acrshort{CPU}s speculate on control flow; could floating-point pipelines analogously speculate on numerical execution and commit results only if accuracy predicates are met?
What is the appropriate speculation window (an \say{error blocking factor}), and what minimal metrics should an augmented \acrshort{FPU} expose to drive such predicates, given that numerical predicates are subtler than branch conditions?

\subsubsection{Co-design metrics beyond toy examples}
With strict hardware budgets, how much historical information about rounding
sensitivity (within a window of operations) is needed to materially improve
algorithm-level stability decisions? Which aggregate metrics are actually useful
and robust across workloads, as opposed to being contrived or overly specific?

\subsubsection{New failure modes in a broad precision landscape}
Low-precision and approximate modes introduce failure modes beyond \acrshort{NaN}s,
overflow/underflow flags, or silent wrong results. What should be signaled, to whom,
and at what abstraction level? Do we need explicit signals for both
\emph{insufficient precision} and \emph{over-computation} (suggesting a safe
downshift)? Can corrective action be performed locally within the \acrshort{FPU}/register
file to avoid expensive checkpoint/restart cycles?

\subsubsection{\acrlong{VaV} under hardware noise}
For approximate computing mechanisms that introduce hardware noise (e.g., voltage
scaling, approximate \acrshort{DRAM}), how should noise models be incorporated into \acrshort{VaV}?
Formal verification exists for floating-point hardware~\cite{Harrison-FP-1999}; can it be extended to certify distributions over computed results and guarantee
desired statistical properties within specified bounds?

\subsubsection{\acrshort{ALU}/\acrshort{FPU} support for statistical validation of \acrshort{QoI}s}
If correctness is increasingly expressed in terms of \gls{QoI} bounds and statistical
validation, can we impose local pointwise error controls that translate into
global \gls{QoI} guarantees? Doing so likely requires cheap verification of candidate
solutions (an optimization problem in itself). Can hybrid correction loops or
surrogate models map local error characteristics to predicted \gls{QoI} bounds
efficiently enough to be used online?

\section{Hardware/Software co-design}
\label{sec:co-design}

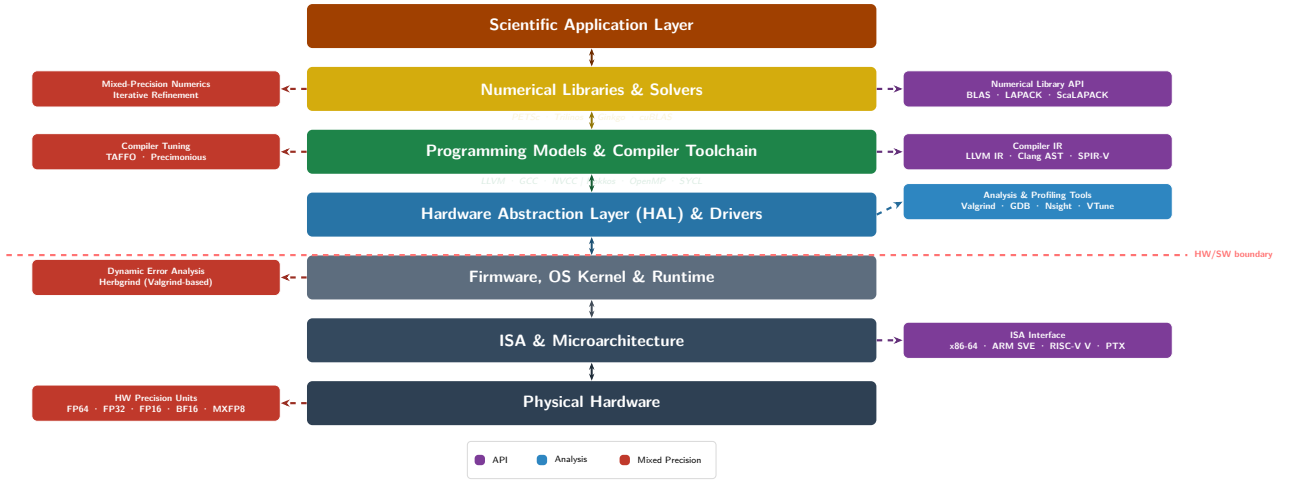
\begin{figure*}[ht!] 
    \centering
    \begin{minipage}{1.0\textwidth}
    \centering
    \resizebox{\textwidth}{!}{
        \begin{tikzpicture}[
          every node/.style={font=\sffamily},
          layer/.style={rounded corners=3pt, minimum width=11cm, minimum height=0.85cm, text centered, text=white, font=\sffamily\bfseries\small, inner sep=4pt},
          iface/.style={rounded corners=3pt, minimum width=5.2cm, minimum height=0.7cm, align=center, fill=ifaceColor, text=white, font=\sffamily\bfseries\tiny, inner sep=3pt, draw=white, line width=0.5pt},
          mpbox/.style={rounded corners=3pt, minimum width=4.8cm, minimum height=0.7cm, align=center, fill=mpColor, text=white, font=\sffamily\bfseries\tiny, inner sep=3pt, draw=white, line width=0.5pt},
          arrow/.style={-{Stealth[length=5pt,width=4pt]}, line width=1.2pt, color=#1},
          dblarrow/.style={{Stealth[length=4pt,width=3pt]}-{Stealth[length=4pt,width=3pt]}, line width=1.1pt}
        ]

        \node[layer, fill=hwColor] (L0) at (0,0) {Physical Hardware};
        \node[layer, fill=isaColor, above=0.35cm of L0] (L1) {ISA \& Microarchitecture};
        \node[layer, fill=rtlColor, above=0.35cm of L1] (L2) {Firmware, OS Kernel \& Runtime};
        \node[layer, fill=drvColor, above=0.35cm of L2] (L3) {Hardware Abstraction Layer (HAL) \& Drivers};
        \node[layer, fill=runtimeColor, above=0.35cm of L3] (L4) {Programming Models \& Compiler Toolchain};
        \node[layer, fill=fwkColor, above=0.35cm of L4] (L5) {Numerical Libraries \& Solvers};
        \node[layer, fill=appColor, above=0.35cm of L5] (L6) {Scientific Application Layer};

        \node[font=\sffamily\itshape\tiny, text=white!85!runtimeColor, below=1pt of L4, yshift=3pt] {LLVM · GCC · NVCC | Kokkos · OpenMP · SYCL};
        \node[font=\sffamily\itshape\tiny, text=white!85!fwkColor, below=1pt of L5, yshift=3pt] {PETSc · Trilinos · Ginkgo · cuBLAS};

        \foreach \a/\b/\c in {L0/L1/hwColor, L1/L2/rtlColor, L2/L3/drvColor, L3/L4/runtimeColor, L4/L5/fwkColor, L5/L6/appColor}
          \draw[dblarrow, color=\c!70!black] (\a.north) -- (\b.south);

        \node[iface, right=0.5cm of L1] (ifISA) {ISA Interface\\x86-64 · ARM SVE · RISC-V V · PTX};
        \node[iface, fill=toolColor, right=0.5cm of L3, yshift=0.25cm] (ifAnal) {Analysis \& Profiling Tools\\Valgrind · GDB · Nsight · VTune};
        \node[iface, right=0.5cm of L4] (ifPM) {Compiler IR\\LLVM IR · Clang AST · SPIR-V};
        \node[iface, right=0.5cm of L5] (ifLib) {Numerical Library API\\BLAS · LAPACK · ScaLAPACK};

        \foreach \src/\dest in {L1/ifISA, L4/ifPM, L5/ifLib}
          \draw[arrow=ifaceColor!80!black, dashed] (\src.east) -- (\dest.west);
        \draw[arrow=toolColor!80!black, dashed] (L3.east) -- (ifAnal.west);

        \node[mpbox, left=0.5cm of L0] (mpHW) {HW Precision Units\\FP64 · FP32 · FP16 · BF16 · MXFP8};
        \node[mpbox, left=0.5cm of L4] (mpComp) {Compiler Tuning\\TAFFO · Precimonious};
        \node[mpbox, left=0.5cm of L2] (mpAnal) {Dynamic Error Analysis\\Herbgrind (Valgrind-based)};
        \node[mpbox, left=0.5cm of L5] (mpLib) {Mixed-Precision Numerics\\Iterative Refinement};

        \foreach \src/\dest in {L0/mpHW, L2/mpAnal, L4/mpComp, L5/mpLib}
          \draw[arrow=mpColor!80!black, dashed] (\src.west) -- (\dest.east);

        \draw[dashed, line width=1.2pt, color=red!55]
          ($(L2.north west)+(-5.8,0)$) -- ($(L2.north east)+(6.0,0)$) node[right, font=\sffamily\tiny\bfseries] {HW/SW boundary};

        \node[draw=gray!30, fill=white, rounded corners=2pt, inner sep=4pt, below=0.3cm of L0] (leg) {
          \begin{tikzpicture}[every node/.style={font=\sffamily\tiny}]
            \fill[ifaceColor] (0,0) rectangle (0.2,0.2); \node[right] at (0.2,0.1) {API};
            \fill[toolColor] (1.2,0) rectangle (1.4,0.2); \node[right] at (1.4,0.1) {Analysis};
            \fill[mpColor] (2.8,0) rectangle (3.0,0.2); \node[right] at (3.0,0.1) {Mixed Precision};
          \end{tikzpicture}
        };

        
        \end{tikzpicture}
    }
    \end{minipage}
    \caption{HW/SW Codesign Stack for Mixed-Precision Computing.}
    \label{fig:codesign_stack}
\end{figure*}

Co-design of a scientific computing system is the iterative, joint design of
hardware, system software, scientific libraries and application frameworks,
and domain science workflows. This process is guided by representative workloads (captured in benchmarks or mini-applications) and measurable targets
such as time to solution, energy, accuracy, and cost.

A central goal is to enable experts at each layer to focus on their own concerns,
while communicating clearly the desires, constraints, and tradeoffs that matter
to other layers. Co-design does \emph{not} imply that any single layer ``gets what
it asks for.'' Rather, it creates a shared vocabulary for what is easier or harder to achieve across the stack and supports informed compromises.

Effective co-design requires a broad range of expertise. Domain scientists, numerical mathematicians, performance engineers, compiler and runtime developers, and computer architects should all participate in jointly optimizing the system.
From the perspective of mixed precision, co-design is the mechanism
by which format choices, architectural features, and numerical guarantees become usable at scale. 

Large-scale scientific computing applications often dictate stringent accuracy and robustness requirements. These rely on basic mathematical algorithms and libraries that employ one or more floating-point formats to provide reliable results regardless of the numerical properties of the input data. Libraries are developed using programming languages and tools that transparently provide access to hardware-level features such as floating-point arithmetic units and instructions, and exceptions. Therefore, as a new feature emerges at the hardware levels, all the layers of this stack (see \cref{fig:codesign_stack}) must adapt to take advantage of it and improve performance without compromising the accuracy or robustness of applications. As a consequence, the objective of co-design is to collaboratively improve all layers of the hardware/software stack to ease development and improve the efficiency of mathematical libraries and applications that use multiple precisions simultaneously. The system should maximize performance and scalability, and
minimize energy, while maintaining guarantees on achievable accuracy and robust
preservation of important invariants.

At the hardware and system-software level, key concerns include communication
intensity, storage hierarchy, and \acrshort{ALU} design. At the programming-system level, we need language and compiler support to expose hardware capabilities, plus
low-level libraries and runtimes to make these capabilities usable. We also want
tooling that can (i) automatically infer or recommend precision in common cases,
(ii) dynamically monitor exceptional conditions to provide guardrails and
debugging support, and (iii) switch from a fast-but-risky path to a more
conservative method when needed. Finally, we need sound mathematical building blocks---backed by error analysis---that provide end-to-end accuracy and
robustness guarantees, together with facilities for experimenting with mixed
precision in application codes to understand domain-relevant accuracy.

\subsection{State of the art}

\subsubsection{Successful precedents in co-design}
Co-design has already succeeded in several scientific computing efforts. While
many systems have historically been designed primarily for peak performance and
cost, there has been increasing emphasis on real workloads represented by full
applications and mini-apps. A particularly successful example is the design of
Japan's Fugaku supercomputer~\cite{Sato2022}, and current efforts include the \acrfull{EPI}\footnote{https://www.european-processor-initiative.eu/codesign-workshop/}.
Outside of scientific computing, modern AI stacks and data centers represent an
evolving co-design effort at enormous scale.

Successful co-design must begin early: silicon is only part of the story, and
a robust software ecosystem is often the determining factor. The co-evolution
of hardware, software, and applications is illustrated by the emergence of
general-purpose \acrshort{GPU} programming. Graphics workloads drove early \acrshort{GPU}s,
but researchers repurposed shader hardware for general matrix computations.
The Brook system at Stanford~\cite{Buck2004} is a notable example: its developer
later contributed to NVIDIA's CUDA. As \acrshort{GPU} acceleration of non-graphics workloads
grew, vendors added hardware support (e.g., improved double precision on some
models) and invested heavily in tuned libraries. In practice, ecosystem maturity
became a major differentiator among competing accelerators. 
Many of these capabilities overlapped with what would later become essential
for modern machine learning workloads.

\subsubsection{Mixed-precision progress enabled by \acrshort{AI}-driven hardware trends}
Motivated by the growth of \acrshort{AI} and machine learning, low and very-low-precision
processing has become widely available on modern \acrshort{CPU}s and \acrshort{GPU}s, while
high-precision performance has stagnated. Scientific computing has nonetheless
leveraged these features---even though they were designed primarily for other
workloads~\cite{hima22}. Successful examples, especially in numerical linear
algebra, include mixed-precision \acrshort{GEMM}, batched kernels with higher-precision
accumulation%
, mixed-precision iterative refinement~\cite{a.b.h.l:24,c.h:18},
and emulation/Ozaki-style techniques for \acrshort{DGEMM}~\cite{10.1145/3773656.3773670}.
These methods are increasingly available in both vendor libraries (e.g.,
cuBLAS~\cite[\S 1.5.2]{cuBLAS}) and academic open-source packages.

\subsubsection{Barriers to broader adoption}
Despite progress, several barriers prevent mixed precision from reaching its
full potential in scientific workloads:
\begin{itemize}
  \item \textbf{Low arithmetic intensity.} Many scientific kernels are
  communication- or memory-bound, so benefits often arise mainly from reduced
  data movement rather than faster arithmetic.
  \item \textbf{Weak exception handling and diagnostics.} Many languages that
  support \acrshort{IEEE} arithmetic lack robust support for floating-point exception
  handling and error detection as specified in \acrshort{IEEE}~754. Even in C, support is
  optional~\cite[\S 7.6]{ISO:2024:III}. At lower precisions, such support becomes
  more critical.
  \item \textbf{Fragility under compiler optimizations.} Predictable floating-point behavior remains difficult when compilers apply algebraic transformations that
  are identities over the reals but not in floating point.
  \item \textbf{Poor functional/performance portability.} Lack of standards and
  rapidly changing architectures make portability difficult. Usable models of
  time and energy costs remain largely the domain of experts and are moving
  targets.
\end{itemize}

\subsection{Future directions}

\subsubsection{Interface-level contracts}

\paragraph*{Error-aware library interfaces.}
To provide flexibility for library designers while relieving application developers from manually selecting formats, we recommend making error characteristics part of library interfaces (see \cref{sec:mp_future_directions}). Designers should specify:
(i) preconditions on input data (e.g., magnitude/range restrictions), and
(ii) accuracy characteristics of outputs. These may include relative/absolute
error bounds, bounds on output magnitudes, or (when appropriate) bounds that
depend on computed quantities such as condition estimators. Different use cases
may require normwise or componentwise specifications, and interfaces should
support both.

\paragraph*{Beyond accuracy: reliability properties.}
Interfaces may also need to describe properties influencing reliable use, such as
reproducibility, preservation of monotonicity or non-negativity, or invariant
preservation. Ideally, users can declare required properties so the system can
select the fastest implementation that satisfies them. This obviously requires developments in multiple components of the hardware/software stack, such as mathematical libraries (see \cref{sec:sw_reproduciibility}), and algorithms (see \cref{sec:emulation}). 

\paragraph*{Compositional checking and optimization.}
A core motivation for interface-level specifications is compositional reasoning.
A compiler or optimization system should be able to verify that a composition of
libraries meets end-to-end error and reliability requirements. For rich interface
languages, checking that one module's guarantees imply another's preconditions
may be nontrivial; demonstrating formal compositional checking is an important
step toward safe, automatic optimization, robust combinations of multiple approximate computing techniques (see \cref{sec:relation}) and use of approximate computing at different levels of hierarchical methods (see \cref{sec:multilevel}). 

\subsubsection{Validation}

\paragraph*{Validation in the absence of full verification.}
Although formal verification of interface properties is desirable, it may not be
feasible universally. A practical intermediate step is to define robust validation
test suites for common interfaces (e.g., \acrshort{BLAS}-like kernels) to check whether
implementations satisfy specified properties.

\paragraph*{Exception-aware execution and guardrails.}
Reliable mixed-precision execution requires detecting exceptional conditions such as invalid operations or overflow/underflow leading to intermediate values outside
normalized ranges. When exceptions can be identified quickly, systems can adopt
a \say{fast and risky} path backed by a slower but more reliable fallback.
Even without full exception handling, surfacing these events is valuable for
debugging and diagnosis~\cite{demmel2022proposedconsistentexceptionhandling}.

Unfortunately, language and compiler support remain incomplete even for mature
\acrshort{IEEE}~754 features such as exception handling and rounding modes. Where languages
support multiple formats, they should also support the associated floating-point
environment, including access to exceptions and rounding modes. At a minimum,
exceptions should be detectable and signaled; ideally, libraries and applications
should degrade gracefully when exceptions occur (see \cref{sec:mp_future_directions}).

Co-design plays a fundamental role in achieving consistent and robust exception handling through research and development in computer architectures (\cref{sec:emerging}), 
floating-point formats (\cref{sec:numberformats}), programming languages and compilers (\cref{sec:sw_lang_comp}, and mathematical libraries (\cref{sec:sw_reproduciibility}).

\paragraph*{Decoupling storage from arithmetic.}
To ease mixed-precision algorithm design and improve memory- and communication-bound
performance, it is beneficial to dissociate storage formats from arithmetic
formats. This requires support across the stack:
\begin{itemize}
  \item \textbf{Hardware:} efficient load/store paths that convert between storage
  types and \acrshort{ALU}-native arithmetic types.
  \item \textbf{Languages/compilers:} ways to specify storage types beyond standard
  \acrshort{IEEE} formats, and to compile type-parameterized algorithms over these types.
  \item \textbf{Runtimes/communication libraries:} awareness of alternate storage
  types for collective operations such as reductions.
\end{itemize}

\paragraph*{Benchmarking for end-to-end accuracy and energy goals.}
Mixed precision is primarily motivated by reductions in time and energy for real
applications. Benchmarks should therefore extend beyond single kernels to include
multi-module workflows with realistic end-to-end accuracy targets. The objective
should not be \say{run a fixed algorithm as fast as possible,} but rather to find
the fastest (or most energy-efficient) configuration that meets prescribed
accuracy goals. Such benchmarks would complement, not replace, micro-benchmarks.
%

\subsection{Open questions}

We have argued for language support for a variety of storage and arithmetic
formats (see \cref{tab:number_formats_desiderata}). Given the large number
of possible formats, co-design across hardware, compilers, and applications is
necessary to determine which subset should be supported and how. Some formats
may be best treated as language-level abstractions but implemented via emulation
rather than native hardware. It is also unclear how to define reasonable default
arithmetic for formats intended primarily for storage.

Energy-aware design is another open area. Energy is not simply \say{joules per
operation times number of operations}: different precision mixes change both
the per-operation cost (and the per-transfer cost; see \cref{fig:Patterson})
and the number/type of operations performed. Hardware telemetry exists, but is
typically coarse-grained and poorly supported by tools. Finer-grained monitoring
could clarify tradeoffs, but measurement alone does not yield predictive models
needed for optimization.

Finally, education and shared understanding are essential. Effective co-design
communication among hardware designers, compiler writers, numerical analysts,
and domain scientists require common high-level models of goals and tradeoffs,
not only low-level format specifications. Developing and disseminating this
shared understanding is inherently cross-disciplinary and requires sustained
support from academia, industry, and policymakers.

\section{Relation to other approximation techniques (precision as one budget among many)}
\label{sec:relation}

Scientific computing (computational science and engineering) workflows combine many approximations. A mathematical model expressed through systems of (in)equalities never fully represents reality; input data (such as material parameters) are often only known with limited precision; discretizations in space and time transform infinite-dimensional models into finite-dimensional ones; and numerical solution and optimization algorithms introduce additional approximations (e.g., low-rank truncation, sparsification, sampling, and tolerance-based termination). Each approximation incurs an \emph{error}, and the dominant error source along the pipeline typically determines the final fidelity of the workflow. In practice, however, individual error contributions are often known only qualitatively (or via pessimistic worst-case bounds), and different stages may quantify error in incompatible ways (different norms, probabilistic vs.\ deterministic guarantees, etc.).

The approximation of the underlying number system by a finite computational
precision (number formats, rounding modes, fixed/floating-point arithmetic)
introduces yet another error source. Historically, applications were executed
in the highest hardware-supported precision, which usually  ``just works''.
As lower and mixed precision become ubiquitous, the interplay between
rounding errors and the many other approximation errors in the workflow
becomes central for obtaining a trusted solution. In this 
section, we discuss current and future directions for balancing precision-induced errors against other
approximation errors at different stages of the workflow.%
%

\subsection{State of the art}
\label{sec:precision_approx_sota}

We briefly review representative recent references that treat finite-precision
rounding errors \emph{together} with other approximation errors. Our aim is
not an exhaustive survey, but to illustrate existing approaches that expose
opportunities for low and mixed precision by balancing the overall error
budget.

\subsubsection{Approximation in the solution process}

\paragraph*{Discretization and multilevel methods}
\label{sec:disc_and_precision}

Work combining discretization errors with finite-precision rounding errors
goes back to Fried's 1971 PhD thesis~\cite{fried1971discretization} and the
1969 report of Melosh and Palacol~\cite{melosh1969manipulation}, although
mixed-precision implementations that \emph{actively balance} these errors were
not considered at the time.

More recently, Tamsdorf et al. \cite{tamstorf_discretization-error-accurate_2021} develop a
progressive-precision scheme for multigrid methods to balance finite-precision
errors with discretization (and algebraic) errors; see
\cref{sec:mp-mg}. Mart{\'\i}nek et al.~\cite{martinek2025exploiting}
introduce a general framework for controlling  computational error in
multilevel sampling methods, treating inexact solves (from early termination
or from mixed-precision iterative refinement) as an additional error component
to balance against discretization and sampling errors. For multilevel Monte
Carlo, they derive conditions and adaptive strategies that set the required
accuracy on each level so that inexactness is safely dominated by the
discretization error.

In time integration, \cite{croci2022mixed} develops mixed-precision
Runge--Kutta--Chebyshev schemes in which the finite-precision error is
maintained at roughly the same order as the time-discretization error. The
schemes use high precision for only a small fraction of function evaluations
without sacrificing overall accuracy.

\paragraph*{Solvers and algorithmic inexactness.}
\label{sec:solvers_and_precision}

In addition to rounding errors, linear solvers are affected by other sources of numerical error and approximation. Some are controlled or prescribed by the
user (e.g., low-rank truncation thresholds in direct solvers, or stopping tolerances in iterative solvers). Others are less controllable and may depend
sensitively on inputs (e.g., unstable pivoting, aggressive dropping in sparse factorizations, or loss of orthogonality in Krylov subspace methods).

A key observation is that many solver-induced inaccuracies can be naturally
\emph{matched} to low or mixed precision: if a solver is intrinsically inexact,
running it in high precision may not yield a commensurate accuracy benefit.
Instead, it can be preferable to run the solver in a lower precision aligned
with its inexactness, thereby improving performance in a numerically harmless
way. This also enables stacking multiple approximation mechanisms that induce
comparable error levels. For example, \cite{abhl23} and \cite{cao_mxpclimate}
combine block low-rank structure, unstable pivoting, and low precision.

For iterative solvers, the same principle applies in terms of accuracy. However, lowering precision can increase iteration counts, reducing or eliminating performance gains. Existing work explores this tradeoff through, e.g., relaxed \acrshort{GMRES} with progressively reduced accuracy in certain operations as orthogonality degrades~\cite{ggl07}; balancing errors from preconditioner applications with the precision at which preconditioners are applied~\cite{cada24,blmv25,bcm25_pre}; combining aggressive restart/tolerances with low precision~\cite{tuwa92,lld22}; balancing the spectral convergence of stationary iterations with subsystem solve precision~\cite{baro15,gvz25}; refining low-precision multigrid solvers (often \acrshort{GPU}-based) to high accuracy
via iterative refinement~\cite{Goeddeke:2007:PAA}; and compressing vectors using tensor-based techniques or black-box compressors~\cite{hal-02572910}.
%

\subsubsection{Approximation in storage and representation}
\label{sec:approximation-in-storage}

%
\label{sec:compression_and_precision}
When a physical signal is sampled at discrete points, then regularity (in a mathematical sense) implies that neighboring samples will have similar values. 
Storing a discretely sampled approximation of a continuous function in a vector will consequently lead to a collection of floating-point numbers that have many redundant bits.
Thus, there is inherently a large potential for using compression techniques.
In Fourier space, the same mathematical features translate into a characteristic asymptotic decay of the Fourier coefficients such that higher frequency modes will contribute less to the signal.
Consequently, higher frequency modes often carry less information that is relevant to the application, and this can offer potential for storing the coefficients in lower precision.
These features are already routinely exploited in compression techniques for image, audio, and video data.
In scientific computing, the use of such techniques is less developed.
Fully exploiting the potential may also require hierarchical representations, see also the discussion in \cref{sec:multilevel}.

Memory accessors separate concerns between data layout in memory and arithmetic on the processing unit. A first class of approaches performs scalar
\lq trimming\rq\ (truncating fraction and/or exponent bits) to reduce memory
footprint and data movement, accelerating memory-bound algorithms. Efficient
implementations have been developed for a range of kernels; see
\cite{gaq21,gjmm24a,ajlm__,kriemann2025hierarchical}.

A second class uses lossy block compression (e.g., ZFP, SZ), originally aimed
at compressing large scientific datasets before post-processing (e.g., in
high-energy physics and planetary missions; see \cite{CAPPELLO2025107323}).
Unlike scalar trimming, these methods exploit intra-block smoothness or
correlation. More recently, they have been integrated \emph{within} numerical
algorithms, including Flexible \acrshort{GMRES}~\cite{hal-02572910}, \acrshort{GMRES}~\cite{hal-03776837},
hierarchical matrix algorithms~\cite{kriemann2025hierarchical}, and sparse
direct solvers~\cite{cgln23}.

A third class comprises problem-structured compression (low-rank, hierarchical,
sparsification), which we discuss next.


\subsubsection{Approximation through structure and randomness}
%
\label{sec:structured_mats}

The growth in scale and complexity of scientific simulations, data analytics,
and machine learning workloads has created intense pressure on storage capacity,
memory bandwidth, and communication. Many applications generate dense or
implicitly dense linear operators (e.g., from \acrshort{PDE} discretizations, kernel
methods, covariance estimation), yet these operators often admit exploitable
structure. Structured and randomized matrix representations approximate large operators while
preserving essential spectral, algebraic or geometric properties, reducing storage,
communication, and arithmetic.

\paragraph*{Globally low-rank matrix approximations.}
Low-rank structure is central when singular values decay rapidly. Low-rank
approximations can be computed via classical \acrshort{SVD} or truncated QR~\cite{GolubVanLoan2013}
and via randomized variants.
In floating-point arithmetic, low-rank structures provide a natural foundation
for mixed precision: dominant singular subspaces may require higher precision,
while less dominant subspaces can be stored/processed at reduced precision, as analyzed
in~\cite{abbg22}.

\paragraph*{Data-sparse hierarchical matrix formats.}
Data-sparse matrix formats generalize low-rank ideas to blockwise low-rank
structure in matrices from Green's functions, boundary integral operators, and
kernel evaluations~\cite{borm_surveyH}. Hierarchical matrices
($\mathcal{H}$/$\mathcal{H}^2$/\acrshort{HODLR}/\acrshort{HSS}/\acrshort{BLR}) enable near-linear storage and
arithmetic complexity~\cite{Borm2003H2,Hackbusch2015,mblr}. These formats create
direct opportunities for mixed precision: near-field blocks can use higher
precision, while far-field low-rank blocks can use reduced/mixed precision in
Schur complements, low-rank updates, basis transfers, and related computations
\cite{abbg22,hong_mxpseismic,cao_mxpclimate,carson2025mixed}. For memory-bound
operations such as matrix--vector products, memory accessors can reduce data
traffic by storing operators in adaptive reduced precisions while computing in
native arithmetic precision (see \cite{ajlm__} for \acrshort{BLR} and
\cite{kriemann2025hierarchical} for hierarchical matrices).

\paragraph*{Sparse matrix factorizations.}
An alternative route to compressing dense matrices is to approximate them as a product of several sparse factors~\cite{gribonval:hal-04954574}.
A particularly important family is \emph{butterfly factorizations}, which combine strong expressivity with an extremely sparse, structured pattern.
Butterfly matrices arise naturally in classical settings---most notably as building blocks in fast transforms such as the Hadamard and Fourier matrices---and they have been leveraged more broadly to accelerate Gaussian elimination (by mitigating the need for pivoting), to design fast structured linear solvers, and to compress neural networks.

In floating-point arithmetic, butterfly factorizations are also well suited to aggressive quantization.
Because they are composed of scaled rank-one updates, they exhibit scaling invariances that can be exploited to quantize factors to very low precision while maintaining a small overall approximation error~\cite{gmr__}.
This makes sparse factorization-based compression a compelling complement to low-rank and hierarchical formats, especially when the primary constraint is memory traffic rather than peak arithmetic throughput.

\paragraph*{Randomized algorithms.}
\label{sec:randomized}

Randomized algorithms are widely used in numerical linear algebra for
dimension reduction~\cite{HalkoMartinssonTropp2011,martinsson2020randomized}.
They employ sketching operators that map large matrices to smaller ones while
preserving target properties with high probability (e.g., approximate norm
preservation of matrix--vector products). Sketching can be done via sampling
(requiring preprocessing) or via random mixing using Gaussian or structured
transforms. Randomized methods underlie least-squares solvers, low-rank
approximations, matrix decompositions, and iterative solvers.

In finite precision, randomization and low/mixed precision exhibit several
synergies. Beyond reducing problem size, randomization can stabilize algorithms
that are otherwise unstable in low precision, e.g., \acrshort{GMRES} with non-orthogonal
bases~\cite{natr24,bcm25}. Adaptive precision strategies are explored for a
randomized range finder in~\cite{chp22} and for truncated QR with randomized
pivoting in~\cite{bmp25}. Sketching in lower-than-working precision (especially
with Gaussian sketches) is employed for randomized \acrshort{SVD}~\cite{ooyo23}, randomized
Nystr\"{o}m~\cite{carson2024single}, and sketched QR decompositions~\cite{georgiou2023mixed,carson2025mixed_sketch}.
Finally, reduced precision can be used throughout randomized \acrshort{SVD} when paired
with higher-precision iterative refinement~\cite{bkmr25,bdkm24}.

\paragraph*{}
Overall, structured representations---including sparse factorizations, hierarchical compression, and low-rank structure and randomized approximations --- form a coherent multiscale compression toolkit.
Combined with mixed-precision arithmetic and storage, these tools are central to scalable solvers and energy-efficient computation on heterogeneous platforms, supporting applications in computational astronomy, climate and weather prediction, seismic imaging~\cite{apps_mxp2022,oaab23}, and genome-wide association studies~\cite{10.1109/SC41406.2024.00012}.
\subsection{Future directions}
\label{sec:precision_approx_future}
In order to minimize the time and energy needed to obtain a trusted solution, we have to
avoid both \lq oversolving\rq\ (wasting resources on
non-dominant error sources) and \lq underresolving\rq\ (incurring unacceptable
errors). Instead, approximation choices---including precision---should be
balanced holistically, end-to-end. Ultimately, this means jointly considering
accuracy targets and resource budgets (time, energy, node count).

We note here that the possibility to use reduced precision is often first a question to the application development and that critical variables are suitably defined.
For example, when weather models would use the Kelvin scale rather than the centigrade scale, almost two decimal digits of accuracy would be wasted with all standard floating point formats.
Similarly if all geometric coordinates were defined relative to the Earth Center, then measuring our typical distances on the meter or kilometer scale would inherently always loose 3-6 digits of accuracy. 
While these examples may sound trivial, we point out that a similar situation occurs when handling large collections of statistical quantities that are all narrowly centered around an expected value.
This is in fact the situation in fluid simulations with the Lattice Boltzmann Method.
Here indeed, a safe use of less than double precision depends foremost on defining and recalibrating all distribution functions around their expected value.
This is, e.g., demonstrated in \cite{Hennig-2023-advanced}. 
Note also that a black-box approach to the problem using memory accessors and compression could reduce the storage overhead automatically, but that this would still not permit the use of lower precision calculations to execute the collision kernels.
This step requires a much more intricate transformation of the algorithm
that can also be automated, see \cite{Hennig-2023-advanced}, but only when exploiting nontrivial application knowledge.

Next, we highlight two broad directions for treating precision as an explicit design variable: improving how we \emph{understand and measure} errors, and improving how we \emph{control} rounding errors to match other error sources better.

\subsubsection{Toward compatible, quantitative error budgets}
\label{sec:error_budgets}

Efficient schemes require a quantitative understanding of error contributions at
each stage of the workflow and, crucially, at interfaces between stages. 
Matching approximation
orders (e.g., spatial vs.\ temporal discretization) does not imply identical
precision requirements for assembly, preconditioning, and solves. 
Especially in hierarchical approximation methods, such as in multigrid methods, a large potential lies in exploiting that not all levels of resolution may require the same numerical precision. 
Especially the incremental contributions of higher resolution, which traditionally dominate storage and computational cost, may be inherently amenable to lower precision representations.

Progress will require tighter a priori and a posteriori bounds for common building blocks,
ideally expressed in compatible norms and incorporating data-dependent
quantities (e.g., conditioning). Automated error analysis and runtime monitoring
could enable tuning and adapting precision during execution.

\subsubsection{A ``precision knob'': better control of rounding error}
\label{sec:precision_knob}

Many approximations admit continuous control (e.g., truncation thresholds in
low-rank methods, stopping tolerances in iterative solvers). By contrast,
rounding error is tied to the unit roundoff of the arithmetic and thus to a
small set of discrete hardware precisions. To balance rounding error against
other error sources and avoid oversolving, we need tools that offer finer
control---ideally a continuum-like \say{precision knob}.  Several promising
routes include:

\paragraph*{More hardware formats.}
The more number formats supported in hardware, the richer the optimization space.
Some formats may be better suited than others; see \cref{sec:numberformats}.

\paragraph*{Adaptively chosen precisions within algorithms.}
Even with a small number of formats, accuracy can be finely controlled by partitioning computations into steps with different sensitivities and assigning precision accordingly. Examples include \acrfull{SpMV}~\cite{doi:10.1137/22M1522619}, low-rank approximations
\cite{abbg22,bmp25}, and tile factorizations~\cite{apps_mxp2022}. Generalizing
this methodology to broader algorithm classes is an important direction.

\paragraph*{Memory accessors for storage-only precision.}
Memory accessors build a practical continuum by using custom storage formats
while decompressing into standard arithmetic formats on the fly. Making such
accessors efficient and widely available remains a key research challenge.

\paragraph*{Emulation for fine-grained arithmetic control.}
Emulation (\cref{sec:emulation}) can provide a relatively fine precision
dial through the number of slices/words used. Emulation-based approaches with
adaptive slice counts appear especially promising for compute-bound settings
involving other approximations, e.g., mixed-precision tile factorizations
\cite{apps_mxp2022} and randomized low-rank approximations~\cite{bmp25}. This may
require emulation designs that decouple multiplication from slicing so that
slicing is performed minimally (ideally once per operand).

\subsubsection{Underexplored opportunities}
\label{sec:underexplored}

\paragraph*{Numerical optimization.}
Despite strong connections to numerical linear algebra and machine learning,
low/mixed precision in optimization remains less explored. Early efforts address
gradient descent~\cite{dfro17} and trust-region methods~\cite{grto20}. Recent work
extends mixed-precision analyses from Newton rootfinding~\cite{ti01,ke22} to
Newton optimization in~\cite{bcmr__}, including interplay with inexact Newton
and Gauss--Newton approximations.

\paragraph*{Tensors.}
Approximate tensor decompositions are common to mitigate exponential growth in
dimension~\cite{bako25}. While truncation-stability is well understood in exact
arithmetic~\cite{bako25}, finite-precision effects are less studied. A general
framework for analyzing errors in tensor operations is proposed in \cite{bkmr__a}.
It is natural to investigate mixed precision for tensors, including extensions
of adaptive-precision low-rank ideas from matrices~\cite{abbg22}.

\paragraph*{Neural networks.}
The propagation of rounding errors during inference and the resulting impacts on accuracy remain incompletely understood.  %
Normwise backward stability for a \acrfull{MLP} under deterministic and probabilistic error models is analyzed in \cite{bbgm__}. 
Authors in \cite{efmr__} develop\st{s} componentwise analyses and propose\st{s} mixed-precision
inference where recomputing a small fraction of inner products in high precision
improves the performance--accuracy tradeoff. Recent work studies transformer
stability for large language models \cite{bfzp25}; extends the \cite{efmr__}
approach to transformers \cite{bgyt26}. Many questions remain, especially when combining
rounding error with other approximations such as sparsification and
quantization~\cite{richatrik_ef21}.

\paragraph*{\acrfull{AMR}.}
Local refinement typically increases condition numbers significantly, suggesting corresponding precision increases; early experiments support this intuition.
Existing rounding-error analyses for \acrshort{AMR} (e.g., \cite{alvarez2012round}) may serve
as a starting point for combining \acrshort{AMR} estimators with precision-aware diagnostics;
see \cref{sec:AMR}, but a full analysis should also consider how hierarchical methods may help to limit the growth of the condition numbers.

\paragraph*{Matrix equations and matrix functions.}
Mixed precision for matrix equations is explored for Lyapunov equations via
iterative refinement~\cite{benner2011mixed,schulze2025towards,benner2025mixed}
and for Sylvester equations~\cite{dmytryshyn2025mixed}. Emerging sketching-based methods for matrix equations suggest additional opportunities for mixed precision.
For matrix functions, arbitrary-precision approaches are studied in \cite{liu2022computing};
here, truncation error (e.g., Taylor approximation) and finite-precision error
should be jointly balanced, and the stopping criteria in iterative schemes raise
similar balancing questions.

\subsection{Open questions}
\label{sec:precision_approx_open}

Beyond the directions above, several challenges in the interplay of precision
with other approximation choices remain largely open. We group them around
(1) how errors are measured, (2) how multiple error sources compose, and
(3) how errors affect algorithmic dynamics.

\subsubsection{Measuring errors: normwise vs.\ componentwise}
\label{sec:measuring_errors}

Many approximation errors (solver, discretization, compression) are bounded
primarily in a normwise fashion, whereas floating-point errors are often
bounded componentwise. Understanding how to combine normwise and componentwise
descriptions---and when one is the appropriate currency---remains an important
gap. Conversely, some floating-point techniques (e.g., emulation) naturally
yield normwise bounds, which may integrate more directly with other error models.

\subsubsection{Composing sources: beyond additive independence}
\label{sec:composing_sources}

A common working assumption is that error sources can be analyzed separately
and combined additively. Future workflows may violate this assumption: the
composition may become order- and context-dependent, and two individually stable
approximations may interact so that their combination is no longer stable under
reduced precision. Identifying such adversarial regimes (and sufficient
conditions that exclude them) is an open problem.

\subsubsection{Analyzing effects on convergence and trajectories}
\label{sec:convergence_trajectories}

For some algorithm classes, attainable accuracy under rounding error is fairly
well characterized, yet the impact on the \emph{convergence trajectory} can be
highly problem-dependent. A canonical example is the \acrfull{CG} method: we can bound attainable accuracy under low/mixed precision, but relating
rounding errors in specific operations to iteration-by-iteration behavior is
notoriously difficult. This matters in practice: if low precision reduces
time-per-iteration but substantially increases iteration count, net gains may
vanish. Developing predictive models for such tradeoffs, and identifying which
operations are truly sensitivity-critical, remains a major challenge.

\paragraph*{}We note that when the precision choice is understood as a component of the total error budget shared among various approximations, expressing contracts for their safe composition becomes an important software challenge.

\section{Software design: user experience, safety, and portability}
\label{sec:softwaredesign}

Software design for reduced- and mixed-precision arithmetic should prioritize user experience, safety, and portability from the start. For users, we need to provide intuitive \acrshort{API}s that abstract precision details, reveal predictable performance and accuracy, and expose clear error diagnosis with automatic fallbacks when needed. For safety, we need to implement numerical safeguards: error bounds, stability checks, consistent rounding behavior, and platform-independent reproducibility. We need to log warnings if precision reductions could alter results beyond acceptable thresholds. For portability, we need to use standard data types (e.g., \acrshort{FP16}, \acrshort{BF16}) and portable mixed precision kernels, ensuring deterministic behavior across architectures, compilers, and \acrshort{GPU} accelerators. Autotuning and robust testing are included, with open benchmarks to verify cross-platform performance. Together, these design choices yield reliable, user-friendly software that harnesses speed and energy savings without sacrificing correctness.

\subsection{State of the art}
\label{sec:soft}

The use of higher precision may be necessary in parts of a computation because
(i) the target solution requires it, (ii) the underlying problem is ill-conditioned,
or (iii) the chosen numerical method is unstable.
However, higher precision is often \emph{not} needed everywhere (and sometimes not even in the final reported solution), creating opportunities to exploit lower precision safely.
Moreover, the accuracy of the input data or the fidelity of the physical model may not justify uniform use of high precision,
so the goal is frequently to match precision choices to the dominant sources of error in the end-to-end workflow. %

\subsubsection{Foundational software building blocks}
\label{sec:mp_building_blocks}

Before reviewing mixed-precision capabilities in numerical libraries, it is useful to outline the software
components that mixed-precision libraries commonly build upon.

\paragraph*{Arbitrary and extended precision arithmetic.}
Representative tools include \acrshort{MPFR}~\cite{mpfr_website}, GNU MP~\cite{gnump_website}, ARPREC~\cite{arprec_github},
double-double and quad-double (DD/QD)~\cite{hida2001,qd_github}, \texttt{chop}~\cite{doi:10.1137/19M1251308},
LoFloat~\cite{lofloat_github}, FloatX~\cite{10.1145/3368086, floatx_github}, and Universal~\cite{omtzigt2023universal}.

\paragraph*{Compression and reduced-storage mechanisms.}
Numerical compression libraries such as ZFP~\cite{zfpwebsite} and SZ~\cite{szwebsite} reduce storage and data movement
and are increasingly integrated \emph{within} numerical algorithms, not only for offline storage.

\paragraph*{Memory accessors.}
Memory accessors decouple \emph{storage precision} from \emph{arithmetic precision} by compressing/trimming data in memory
and converting on the fly at load/store boundaries, see also \Cref{sec:approximation-in-storage}.
Examples include Ginkgo~\cite{https://doi.org/10.1002/spe.3041}, a \acrshort{BLAS}-based block memory accessor in MUMPS~\cite{ajlm__},
and related efforts.

\subsubsection{Library support for multiple and mixed precisions}
\label{sec:mp_library_landscape}

%
\paragraph*{A fixed menu of precisions.}
Historically, numerical libraries have frequently offered a small, fixed set of hard-coded precisions.
This gives users explicit control, but maintaining multiple versions is error-prone, and implementations can drift out of sync.
Moreover, guaranteeing execution in a given precision does not guarantee a given \emph{solution accuracy}.

Examples providing \acrshort{FP32} and \acrshort{FP64} variants include
\acrshort{BLAS}~\cite{10.1145/355841.355847,10.1145/42288.42292,10.1145/77626.79170},
LAPACK~\cite{lapack_users_guide:1999}, PSCTOOLKIT~\cite{psctoolkit},
QR-MUMPS~\cite{10.1007/978-3-642-40047-6_53}, MUMPS~\cite{10.1007/3-540-70734-4_16,mumps_website},
HSL~\cite{hsl}, SuperLU~\cite{doi:10.1137/S0895479895291765},
PastiX~\cite{10.1016/S0167-8191-01-00141-7}, and NAG~\cite{nag}.

\paragraph*{Templating and precision-generic libraries.}
Templated libraries aim to make algorithms \say{future-proof}: write the algorithm once and instantiate it at multiple precisions.
Examples include
\cpp~mdspan~\cite{Hollman_2019}, Kokkos~\cite{9973611}, \tlapack~\cite{tlapack},
deal.II~\cite{10.1145/1268776.1268779}, SLATE~\cite{gates2020slate}, STRUMPACK~\cite{10.1145/2930660},
Ginkgo~\cite{ginkgo_website}, Composyx~\cite{composyx_gitlab}, and HyTeG~\cite{kohl2019hyteg}.
Fortran supports multiple precisions through \texttt{kind} parameters and can be used to write precision-agnostic code,
even though it lacks C++-style templates; in practice, compiler support for nonstandard or emerging data types may still limit portability.
In C++23, half precision is part of the standard (\texttt{std::float16\_t}).

\emph{Limits of templating in large production codes.}
Experience from large geoscience modeling efforts suggests that templating alone may be insufficient for \acrshort{GPU} targeting and mixed-precision execution.
In the ICON model and the HyTeG matrix-free multigrid framework~\cite{kohl2019hyteg, bohm2025code},
the complexity of refactoring a flexible, general-purpose infrastructure to support \acrshort{GPU}s and multiple/mixed precisions
led to a shift in strategy: development moved toward a specialized application focused on Earth-mantle geometry and physics.
This illustrates a recurring tension: generality and extensibility can become too costly when the primary constraints are \acrshort{GPU} portability
and rapid mixed-precision enablement under limited development resources.

\paragraph*{Mixed precision: growing, but often ad hoc.}
Mixed precision has entered scientific software in a limited, sometimes ad hoc manner, typically through targeted optimizations.
Libraries with some mixed-precision capabilities include
BLIS~\cite{10.1145/3402225}, LAPACK~\cite{lapack_users_guide:1999}, cuBLAS~\cite{10.1145/3773656.3773670,nvidia:Vishwanath:2019},
XBLAS~\cite{xblas}, Ginkgo~\cite{ginkgo_website}, SuperLU~\cite{doi:10.1137/S0895479895291765},
STRUMPACK~\cite{10.1145/2930660}, MUMPS~\cite{10.1007/3-540-70734-4_16,mumps_website},
Composyx~\cite{composyx_gitlab}, AdaptSPMV~\cite{doi:10.1137/22M1522619},
PLASMA~\cite{plasma-users-guide-2009},
MAGMA~\cite{tdb10,tnld10,dghklty14,dghklty14}, {\tt HSL\_MA79}~\cite{hogg2011article},
HPL-MXP~\cite{hpl-xmp}, PETSc~\cite{petsc}, Trilinos~\cite{trilinos}, and HiCMA~\cite{10.1109/SC41406.2024.00012}.
Many Gordon Bell Prize-winning codes employ mixed precision.

For deep learning, mixed precision is routine; \acrshort{CPU} examples include
FBGEMM~\cite{FacebookGEMM}, XNNPACK~\cite{XNNPACK}, and QNNPACK~\cite{QNNPACK}, with \acrshort{GPU} analogs available.

\paragraph*{Application evidence: precision need not be uniform.}
Application codes have traditionally run in a fixed precision chosen during software development — historically, \acrshort{FP64} by default.
Yet many workflows do not require \acrshort{FP64} end-to-end.
For example, an astronomy pipeline at Paris Observatory ran in \acrshort{FP64} by default even though an \acrshort{FP32} implementation was available. Subsequent reviews
found that \acrshort{FP32} is sufficient because the pipeline ultimately drives actuators whose input format is unsigned 16-bit, making \acrshort{FP64} (and even \acrshort{FP32} everywhere) unnecessary~\cite{8945098}.

A seismic imaging application provides a complementary illustration of \emph{spatially varying} precision needs.
A tile-centric mixed-precision strategy reduces memory pressure while preserving accuracy~\cite{eage:/content/papers/10.3997/2214-4609.2025643027}: most tiles can be computed in 16-bit precision without violating the error budget for the final image,
while tiles near the propagating wavefield must remain in 32-bit precision.
The tile-precision selection follows the methodology in Section 14.3 of Higham and Mary~\cite{hima22},
enabling combinations of \acrshort{FP64}/\acrshort{FP32}/\acrshort{FP16} and even \acrshort{FP8} on Grace Hopper hardware.
Vendor \acrshort{BLAS}/LAPACK kernels perform the numerical work, and PARSEC performs on-the-fly type conversions (sender- or receiver-side)
to reduce communication overhead.
This experience highlights a broad challenge: automatically identifying regions of a computation that can be safely demoted.

Other application packages that employ mixed precision include
GROMACS~\cite{gromacs2015} (molecular dynamics),
WRF~\cite{wrf2008} (weather/climate),
NAMD~\cite{namd2005} (biomolecular simulation),
QMCPACK~\cite{qmcpack2020} (quantum Monte Carlo),
LAMMPS~\cite{lammps2022} (materials science),
HPGMG~\cite{hpgmg2014} (PDE multigrid),
OpenFOAM~\cite{openfoam} (CFD),
SPECFEM3D~\cite{specfem3d2002} (seismology),
and AlphaFold~\cite{alphafold2021} (ML for science).

\subsubsection{Libraries using memory accessors.}
Examples include MUMPS~\cite{mumps_website}, Ginkgo~\cite{ginkgo_website},
Composyx~\cite{composyx_gitlab, hal-02572910, hal-03776837}, and HLIBpro~\cite{hlibpro_website}.

\subsubsection{Tools for precision exploration and automated tuning}
\label{sec:mp_tools_precision_exploration}

A key research avenue is the design of tools that, given a desired solution accuracy,
automatically identify where reduced precision is safe.
An example is PROMISE~\cite{chen:hal-05291790}, which provides task-specific validation for automated precision tuning.
From a developer's perspective, the PROMISE workflow can be summarized as: (i) instrument the code to tag variables whose precision may be modified and to specify results to be checked;
(ii) given an acceptable tolerance, use a dichotomic search (Delta Debugging) to determine which variables can be computed and stored in reduced precision,
potentially using custom formats; and
(iii) validate proposed changes via numerical sensitivity analysis.
Internally, PROMISE relies on CADNA~\cite{Eberhart_REC_2015} and has been extended to customized precision through FloatX~\cite{floatx_github}.

Existing tools such as Precimonious~\cite{precimonious-sc13,precimonious_web}, Verrou~\cite{verrou_web},
PROMISE~\cite{promise_web}, Verificarlo~\cite{verificarlo_web},
EXCVATE~\cite{excvate-arith2025}, FPChecker~\cite{fpchecker_web}, and RAPTOR~\cite{raptor_web} provide valuable building blocks, but a unified, user-friendly pipeline that takes users from exploration to deployment remains a major gap.
Another important direction is tooling that supports automatic emulation of a requested precision on hardware that does not natively provide it.
\Cref{tab:tools_capabilities} surveys the important features provided by several tools, focusing on usability from the users' viewpoint.
Most of the tools need LLVM/IR bitcode or source code, hence, they cannot support closed-source executables, such as vendors' mathematical libraries for which we often only have access to the binaries.
The parallel programming support with \acrshort{MPI} and \acrshort{GPU} is still very limited. 

\begin{table*}[h!] 
\centering
\caption{Characterization of the tools for floating-point precision analysis and automated tuning.}
\label{tab:tools_capabilities}
\scriptsize
\bgroup
\renewcommand{\arraystretch}{1.25}
\begin{tabular}{|p{2.0cm}|
p{1.5cm} |
p{1.8cm}|p{1.8cm}|
p{1.5cm}|
p{1.5cm}|
p{2.0cm}|
p{1.5cm}|}
  \hline
&\textbf{EXCVATE}   &\textbf{FPChecker} &\textbf{Precimonious} &\textbf{PROMISE} &\textbf{RAPTOR} & \textbf{Verificarlo} &\textbf{Verrou} \\ \hline
functionality & Spoofs FP exceptions, finds bugs 
& Detects FP exceptions
& Auto-tunes FP variables precision
& Auto-tunes FP variables precision
& Profiles FP precision impact
& Checks FP accuracy stochastically
& Randomizes rounding \& perturb, finds instabilities \\ \hline
compiler based?      &\no, binary instr.         & clang++     & LLVM/IR   & need source code   &  LLVM/IR         & LLVM/IR              &\no, binary instr. \\ \hline
internal arithmetic     & cvc5 \acrshort{SMT} solver  & non     & long double    & \acrshort{DSA}             & \acrshort{MPFR}, native HW types           & \acrshort{MCA}, \acrshort{MPFR}  & \acrshort{CESTAC}, \acrshort{MCA} \\ \hline
language support         & language-agnostic     & C/C++            & C/C++   & C/C++, Fortran  &  LLVM-supported languages   & LLVM-supported languages  & LLVM-supported languages\\ \hline
MPI support?    &  \no         & \no       & \no             & \yes     &  partial    & \no                  & \yes \\ \hline
GPU support?   &  \no     & NVIDIA         & \no                  & partial         &   AMD, NVIDIA   & \no                  & \no \\ \hline
\end{tabular}
\egroup
\end{table*}

\subsection{Future directions}
\label{sec:mp_future_directions}

Mixed-precision computing is rapidly maturing, but the scientific software ecosystem still lacks the abstractions,
tooling, and runtime intelligence needed to make it routine.
The needs expressed by domain scientists, library developers, compiler experts, and hardware vendors suggest several
directions, summarized below.

\subsubsection{Accuracy-aware interfaces that capture user intent}
Users often know \emph{what accuracy they need}, not \emph{what precision to request}, and the proliferation of formats
makes precision-selection an unrealistic burden.
Interfaces should therefore allow users to specify their intent:
desired accuracy (normwise and/or componentwise), acceptable tolerances and stopping criteria,
and priorities such as performance, energy, or reproducibility.
Separately, users can reasonably be asked to provide input/output ranges, which can support range-aware safeguards (e.g., blocking values within exponent bounds).
On the developer side, such intent can map to mechanisms such as \acrshort{BLAS}-like interfaces with an $\varepsilon$ parameter,
memory accessors that compress/decompress automatically (especially important for memory-bound phases),
and hardware-assisted precision selection (especially relevant for compute-bound kernels).
These abstractions support separation of concerns: users specify intent, while libraries/compilers determine how to meet it.
They also motivate testing and \say{grade} implementations against declared contracts.

\subsubsection{Clarifying roles across the software stack}
Mixed precision spans many layers of the software stack: applications, domain frameworks, libraries (\acrshort{BLAS} and beyond), programming models,
compilers/runtime systems, and hardware.
Progress depends on clear boundaries and communication mechanisms that translate requirements across layers without overwhelming users.
In particular, the \acrshort{BLAS}/kernel layer remains a central locus of activity, tooling, and performance engineering,
while co-design at the compiler/runtime and architecture levels determines what precision features can be exposed cleanly and portably.

\subsubsection{Language and compiler support for precision as a first-class concept} \label{sec:sw_lang_comp}
Fortran and C++ offer constructs for controlling floating-point behavior (e.g.,~\cite{cppisnanwebsite}), but support remains inconsistent across compilers and targets.
Mixed-precision software increasingly needs portable, fine-grained control over rounding modes, exception flags,
conversions, and arithmetic formats, ideally without global side effects.
Support for emerging formats (\acrshort{FP8}, \acrshort{BF16}, \acrshort{TF32}, quad precision) and for emulation of non-native formats will be essential.
This represents a shift from historical library design choices (e.g., LAPACK avoiding advanced language features for portability): richer language and compiler support are now needed to express accuracy requirements, error-handling strategies, and precision contracts.
However, it remains a formidable challenge to provide adequate language and compiler support for the vast number of floating point formats surveyed in~\Cref{sec:numberformats},

\subsubsection{``Responsibly reckless'' algorithms and runtime adaptivity}
A recurring algorithmic pattern is \say{try low precision first}, paired with robust failure detection, escalation, and certification, followed by the necessary higher-precision computation in certain parts.
Realizing this reliably requires runtime signals and policies, such as:
residual-based triggers, stagnation indicators, loss of orthogonality measures,
adaptive inner- and outer-iteration schemes and selective recomputation when low-precision steps degrade accuracy.
Such mechanisms are already visible in incomplete factorizations, iterative refinement, and Krylov subspace methods;
extending them to iterative eigensolvers (e.g., Lanczos- or Arnoldi-based), least-squares solvers, and nonlinear methods is a promising direction.
More broadly, progress will require co-designed support for emulation frameworks, domain-aware refinement (analogous to \acrshort{AMR}),
standardized metadata/contracts for routines, and hardware support for mechanisms such as memory accessors and stochastic rounding. %

\subsubsection{Reproducibility, debugging, and exception handling}
\label{sec:sw_reproduciibility}
Reproducibility in mixed-precision computing is best treated as a graded, multiclass service rather than a binary property.
A useful spectrum includes:
(i) bitwise reproducibility (highest cost),
(ii) numerical reproducibility (agreement of key quantities of interest within bounds),
and (iii) statistical reproducibility (distributional agreement for ensemble/stochastic workflows).
Supporting this spectrum requires exception-aware execution models in which hardware flags provide rapid feedback
when assumptions of the underlying error analyses are violated.
Long-term reproducibility also benefits from environment capture tools such as Nix~\cite{nixwebsite} and Guix~\cite{guixwebsite}.
Mixed-precision routines should expose explicit contracts describing guarantees (accuracy, convergence behavior, reproducibility class)
and assumptions (conditioning, scaling, valid input ranges).
Finally, debugging and profiling tools must evolve to reveal precision-related diagnostics, exception tracking,
and interactions between arithmetic formats and algorithmic stability.

\subsubsection{Tools for analyzing use of energy}
\label{sec:sw_energy}
It is increasingly recognized that an important environmental issue for any software is the energy it consumes while running. Some tools are becoming available to help analyze energy usage. An example in cloud computing is RAPL (Running Average Power Limit) from Intel~\cite{rapl}. A review covering tools for carbon footprint estimation, including online calculators, embedded packages, and server-side tools, can be found
in \cite{carbon-footprint}.

\subsection{Open questions}
\label{sec:mp_open_questions}

Several foundational questions remain on the path toward making mixed precision a safe and portable part of scientific software. The first challenge is communication across the layers in the codesign stack in~\cref{fig:codesign_stack}, as decisions regarding precision span applications, libraries, compilers, runtimes, and hardware. Where should the abstraction boundaries lie so that domain scientists can specify intent, while library developers and hardware vendors handle representation choices? We believe some degree of separation of concerns is necessary, as users, developers, and vendors cannot be expected to track the expanding set of formats and representations across this stack, but must be able to express expectations and requirements across layers. In this context, determining a minimal set of assumptions across formats and what guarantees can be offered is crucial, with particular care required for formats lacking exception handling. Thus, one open question is how libraries can expose certified mixed-precision routines while keeping interfaces comprehensible and scalable. 

Beyond interfaces, there is a gap between the successes of mixed- and reduced-precision libraries and the challenges of testing and evaluating these approaches in large applications. How far can code-generation approaches help prototype and evaluate mixed-precision alternatives for such cases? Can the impact of mixed precision be assessed on large applications without relying on exhaustive case-by-case testing?

Finally, to encourage adoption and enable fair comparisons, a set of minimal, shared reference implementations should be standardized, and the tooling environment must evolve to meet the needs of mixed-precision computing. We envision a tool that encompasses everything from analysis, validation, autotuning, profiling, and debugging to the validation and reporting of energy savings, without compromising correctness or interpretability.

\section{Precision as a multilevel resource}
\label{sec:multilevel}

Hierarchical structure is a defining feature of many efficient algorithms and data representations in scientific 
computing, such as \acrfull{MG}~\cite{trottenberg_multigrid_2001}, \acrshort{AMR}~\cite{berger-amr-84}, hierarchical bases~\cite{yserentant1986multi}, hierarchical matrices~\cite{hackbusch_hierarchical_2015}, or wavelets~\cite{stollnitz_wavelets_1996}. These algorithms are the workhorse for many efficient implementations, as discussed in \cref{sec:relation}.
This hierarchy provides a natural coordinate system for assigning
precision. With data and computation being organized across
a multilevel hierarchy, these same structures can be used
to index where and how numerical accuracy is expended. Coarse
levels, fine scales, localized features, and compressed
representations each occupy distinct positions in the
hierarchy, and therefore admit distinct precision regimes
aligned with their algorithmic role.

In this view, precision is no longer a single global choice
or a localized kernel-level optimization. Instead, it becomes
a resource that is distributed across levels, regions, and
representations. Decisions about storage format and arithmetic
precision is made in relation to the, possibly adaptively evolving, 
hierarchical structure of the method, allowing accuracy to be allocated
where it has the most impact or is required, and relaxed where redundancy 
is present.

In this section, we cover case studies of three categories of numerical methods that exploit precision as a multilevel resource:
multigrid methods, adaptive mesh refinement techniques, and hierarchical bases.

\subsection{State of the art}

\subsubsection{Multigrid}

\begin{figure*}
    \centering
    \includegraphics[width=.7\linewidth]{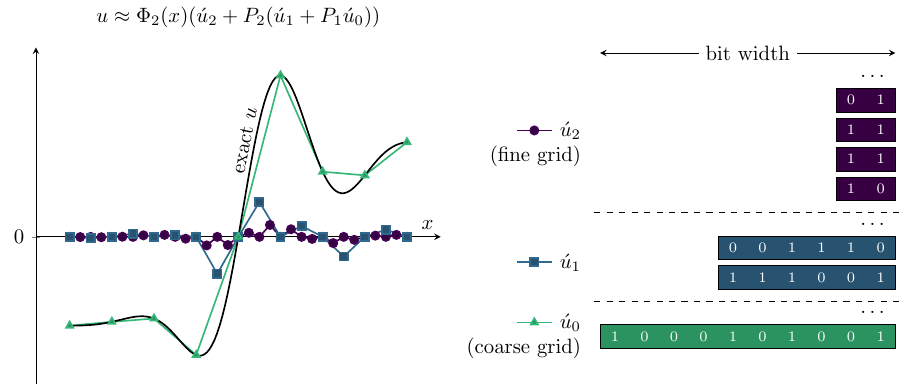}
    \caption{In compact multigrid, the exact solution function $u$ is approximated by a coarse approximation $\acute u_0$ and successively finer corrections (here: $\acute u_1$, $\acute u_2$). After prolongating all sections to the finest grid via prolongation operators $P_1$ and $P_2$, summing, and multiplying with the basis row-vector $\Phi_2$, the exact function is approximately recovered. Due to their smoothness, the corrections decrease in magnitude, allowing them to be stored at lower bit widths on finer grids without loss of accuracy.}
    \label{fig:cmg}
\end{figure*}

When solving \acrshort{PDE}s numerically, the condition number of the discrete differential operator typically increases under refinement.
As a result, fine grid quantities and computations are more susceptible to quantization and rounding errors than coarse ones.
This naturally suggests a vertical precision hierarchy in which coarse levels are computed at lower precision, while finer levels use high precision.
Theoretical analysis~\cite{mccormick_algebraic_2021,tamstorf_discretization-error-accurate_2021} explains the precise dependence of the required precision on the level.
Notably, the precisions of the key \acrshort{MG} quantities (system matrix, solution vector, right-hand side, residual) increase with different rates toward finer levels, meaning that every step of the algorithm effectively has its own precision requirement.
In practice, this results in performance benefits from using low precision on coarse grids~\cite{buttari_block_2022,tsai_three_2023}.
Mixed-precision multigrid with user-specified level-wise precision is implemented for instance in the Ginkgo~\cite{ginkgo_website}, hypre/BoomerAMG~\cite{henson_boomeramg_2002}, or Trilinos~\cite{trilinos} libraries. 
For an overview of mixed precision in software design, see~\cref{sec:softwaredesign}.

While the above approach allows savings on the coarse grids, a large portion of memory and compute requirements stems from the large size of the finest grid.
In addition to storing an increasing number of degrees of freedom, achieving discretization error requires more bits per degree of freedom as the grid is refined.
At the extreme scale, this superlinear growth in memory volume challenges the memory and bandwidth constraints of emerging hardware (cf.~\cref{sec:emerging}), necessitating an efficient, hierarchical representation of the solution vector. 
To mitigate the growth in bit width,  \acrshort{CMG}~\cite{pan_compact_1992,thurner_ganzzahlige_1994,bauer_multigrid_2025} leverages the smoothness of the solution to compress its representation in memory, see also
\cref{sec:approximation-in-storage}.
Instead of viewing multigrid levels as independent representations of the solution, the solution vector is decomposed into a coarsest-grid approximation and successive corrections on the fine grids, see \cref{fig:cmg}. 
This is similar to how fine grids serve to improve the solution that is computed on coarser grids in the Full Approximation Scheme formulation of \acrshort{MG}~\cite{brandt_multigrid_2011}.
Under suitable assumptions on the regularity of the solution, the finest-grid correction can be represented with a constant and surprisingly small number of bits, independent of the total number of levels.
As a result, the full compressed solution can be stored using a constant number of bits per degree of freedom, i.e., with linear storage complexity in bits.
Nevertheless, the precision used in the compact format typically varies from one level to the next (growing toward coarser grids).

Beyond level-wise precision assignment, precision can be further distributed among individual multigrid components~\cite{vacek_mixed_2025}. Different parts of the algorithm (smoothers, residual evaluations, restriction and prolongation operators, and coarse-grid solvers) exhibit different sensitivity to rounding errors and may therefore benefit from distinct precision choices. 

Finally, \acrshort{MG} can be endowed with a second, orthogonal hierarchy, by choosing a number representation with shared or hierarchical exponents, cf.~\cref{sec:numberformats}.
When data is stored in very low precision, these representations avoid excessive storage cost associated with exponents, while giving potential to leverage modern specialized hardware.
In particular, \acrfull{BFP} has been proven suitable~\cite{kohl_multigrid_2024,sundriyal_adaptive_2025}.

\subsubsection{\acrlong{AMR}}\label{sec:AMR}

\acrshort{AMR} dynamically adjusts computational resolution to concentrate degrees of freedom in regions where the solution exhibits sharp gradients, singularities, or other localized features. Originating in the finite difference and finite element communities~\cite{berger-amr-84,berger1989local}, \acrshort{AMR} improves numerical efficiency by either reducing the number of elements required for a target accuracy or minimizing error for a fixed computational budget. Central to \acrshort{AMR} are error indicators or estimators, which guide refinement and coarsening. These include residual-based estimators~\cite{VERFURTH199467}, recovery-based methods such as the Zienkiewicz–Zhu estimator~\cite{zienkiewicz1992superconvergent}, hierarchical estimators based on enriched spaces, and the most recent patch-based mesh refinement in space and time~\cite{berger-amr-24}. 
All provide spatially localized information about the discretization error.

This mesh adaptivity naturally suggests a connection to a hierarchy of precision. Regions flagged for refinement, whether via smaller mesh size (h), higher polynomial degree (p), or locally adapted mesh relocation (r), often demand greater numerical accuracy in operator evaluation and iterative solution. 
In principle, the information from indicators or error estimators could also be used to guide the choice of precision in different parts of the algebraic solver. 

Despite extensive work on \acrshort{AMR} and, separately, on mixed-precision solvers and hierarchical arithmetic, there is currently limited work that combines \acrshort{AMR} with spatially varying precision, although e.g.~the hierarchical strategy of \cite{Ruede-1993-Fully-Adaptive} could be naturally extended to include variable floating point precision.
Existing mixed-precision approaches typically assume uniform precision per solver phase or kernel (not counting recent efforts on adaptive precision), rather than precision distributions driven by mesh-refinement patterns. At the same time, \acrshort{AMR}’s dynamic memory layout has limited adoption on \acrshort{GPU}s, although heterogeneous architectures now allow asynchronous execution of adaptivity on \acrshort{CPU}s while solvers run on accelerators. Recent work on \acrshort{CPU}–\acrshort{GPU} scheduling for block-structured Cartesian \acrshort{AMR} (e.g., AMReX~\cite{Zhang21AMReX}) shows that, without careful co-design, \acrshort{AMR} overhead can limit \acrshort{GPU} throughput; however, asynchronous off-loading and overlapping computation with mesh management enable performance gains when properly integrated. Together with the emerging \acrshort{GPU}-enabled solver stack, this suggests that integrating \acrshort{AMR} indicators with the selection of hierarchies of precisions in such settings represents a largely unexplored and promising research direction.

\subsubsection{Hierarchical bases}

The hierarchical basis method is particularly well-conditioned. Bank et. al. \cite{10.1007/BF01462238} showed that the method is
${\mathcal O}(\log(n))^2$ for problems defined in a
two-dimensional domain. In contrast,
Ong's~\cite{doi:10.1137/S1064827594276539} analysis
suggests that the method is ${\mathcal O}(\log(n)n)$
for problems defined in three-dimensional domains.
Owing to this structure, a low-precision preconditioner combined with a Krylov subspace method is expected to perform robustly: the preconditioner can operate in reduced precision, while the outer Krylov iterations---such as matrix–vector products---are executed in higher precision.
Ruda et al. \cite{doi:10.1177/10943420221084657} demonstrated that this strategy is effective when a direct solver is used as the preconditioner in lower precision. At the same time, most of the precision-sensitive work can remain in the outer Krylov iteration.

\subsection{Connections to other topics}

In the context of \acrshort{AMR}, increasing the polynomial degree in hp-adaptivity often amplifies conditioning effects and necessitates higher numerical precision. Therefore, precision selection should be directly coupled to local polynomial order and refinement patterns. \acrshort{BFP} represents another plausible approach, particularly for max-norm control within blocks, since it uses a shared exponent; tools from fixed-point analysis may provide a theoretical foundation for studying its stability and convergence properties.

The choice of norm is also central in precision analysis. While \acrshort{CG} is naturally analyzed in the energy norm and \acrshort{GMRES} typically in the $\ell_2$-norm, the energy norm often exhibits smoother error behavior and more structured Galerkin projections, potentially simplifying the analysis. Moreover, optimal convergence in the energy norm does not necessarily translate to $\ell_2$, yet it may allow greater tolerance for reduced precision.

In the \acrshort{AMR} setting, fundamental questions arise regarding global matrix storage and whether matrix blocks can be stored with heterogeneous precisions. \acrfull{SEM}, being matrix-free, are particularly promising, since refinement typically affects only selected domains, allowing localized adjustment of polynomial degree and precision. 

\subsection{Future directions}

\subsubsection{\acrfull{CMG}}\label{sec:mp-mg}

\acrshort{CMG}~\cite{pan_compact_1992,bauer_multigrid_2025} provides a promising framework for using reduced precision even on the fine, expensive grids in the hierarchy.
Several research directions appear particularly important for the coming years.

First, it remains to be determined under which conditions the \acrshort{PDE} solution can be represented compactly with a constant bit count.
Similarly, theory that explains the required precision for other quantities, such as residuals and corrections, and the effect of reduced precision on convergence rates, is largely unexplored.

Another challenge is identifying suitable subspaces that enable a clean separation between coarse and fine components.
Ideally, the space of fine-grid corrections should be complementary to the coarse approximation to avoid redundancy in the representation.
A promising candidate are wavelet-based multiresolution basis constructions~\cite{daubechies_ten_1992,stollnitz_wavelets_1996} although it is unclear which type of wavelets integrates best with (compact) \acrshort{MG}.
Some open questions are whether subspaces corresponding to the levels should be mutually orthogonal, whether prolongation matrices should have orthogonal columns, and whether orthogonality should be measured in $L^2$ or in energy~\cite{sudarshan_operator-adapted_2005}.
In any case, there are other options beyond wavelets that achieve (some of) these desirable properties.
For example, the prolongation matrices in the hierarchical-matrix \acrshort{MG} method from~\cite{sushnikova_h2-mg_2025} have orthogonal columns.

Apart from these theoretical and algorithmic improvements, the realization of proof-of-performance implementations is important.
Demonstrating that \acrshort{CMG} can be implemented efficiently on modern, massively parallel architectures without excessive synchronization or communication will be critical to its practical adoption.
This includes assessing whether the algorithmic gains from reduced precision and storage translate into end-to-end performance improvements at scale.

Finally, \acrshort{CMG} currently combines low-precision storage with high-precision computation. A natural next step is to explore how low-precision arithmetic can be employed throughout the \acrshort{CMG} cycle, including residual calculation and smoothing, while preserving discretization-error accuracy and convergence.
Closely related are low precision stencil/matrix representations and the elimination of small (number of elements) yet high-precision temporary vectors, further reducing data movement.

\subsubsection{\acrshort{AMR}}
This is a very appealing and largely unexplored direction.
A natural next step is to couple hierarchies of precision, e.g., in the representation of basis functions and operators, with classical mesh adaptivity strategies such as r-refinement (mesh redistribution) and hp-refinement (simultaneous variation of the mesh size h and the polynomial order p).

Regions subject to strong refinement, particularly under hp-adaptivity, typically demand higher numerical accuracy due to increased polynomial degree, sharper gradients, or localized stiffness. However, this need for higher precision is typically local rather than global. Rather than executing the entire solver stack in uniform high precision, one can propagate mesh adaptivity information to the algebraic solver levels. This could enable, for example, the design of a domain decomposition-based preconditioner in which the precision used in subdomain solves varies across the computational domain, as proposed in \cite{cjmn__}, with refinement indicators guiding the choice of local precision.

Such a strategy naturally leads to solvers with a hierarchy of precisions, in which floating-point formats are adapted in concert with \acrshort{AMR}. Since \acrshort{AMR} usually affects only a relatively small portion of the mesh, its impact can be localized within matrix-based, explicit, or matrix-free solvers without enforcing global changes. In particular, matrix-free formulations, especially in spectral element methods (SEM), offer strong potential for integrating hierarchies of precision directly at the operator evaluation level, where both arithmetic precision and storage formats can reflect local refinement patterns.

\subsubsection{Hierarchical bases}

An advantage of hierarchical basis methods, emphasized by Zienkiewicz et al. \cite{ZIENKIEWICZ198353}, is that the construction of the basis functions inherently provides an error indicator. This indicator can guide both mesh refinement and precision selection by identifying regions where higher precision is needed to capture changes in the solution. In this way, error estimation becomes a unifying mechanism linking discretization adaptivity and precision adaptivity.

Further work is required to realize mixed-precision strategies in combination with hierarchical bases fully. The iterative solver proposed by Bank et al. \cite{10.1007/BF01462238}
resembles an \acrshort{MG} method, but the smoother is applied only to a subset of nodes at each level. Consequently, not all grid levels exchange coarse-grid information in the standard way, and classical \acrshort{MG} precision-transfer strategies cannot be directly reused. An alternative is Griebel's generating systems method \cite{doi:10.1137/0915036}, which exhibits similarly favorable conditioning in both two- and three-dimensional settings and whose \acrshort{MG} equivalent more closely resembles standard \acrshort{MG} algorithms. This makes it a promising framework for systematically embedding hierarchies of precision into adaptive solvers.

\subsection{Open questions}

Treating numerical precision as a multilevel, spatially and algorithmically varying resource raises several fundamental open questions that extend beyond any single class of solvers.

A first and overarching question concerns the magnitude of the achievable gains. While hierarchies of precisions promise reductions in memory footprint, bandwidth, and arithmetic cost, it remains unclear how much of this potential can be realized in end-to-end applications. Quantifying these gains and identifying when they are asymptotically meaningful rather than merely constant-factor improvements remain open problems across hierarchical numerical methods.

Closely tied to this is the role of hardware (see \cref{sec:emerging}). It is not yet clear whether the benefits of precision hierarchy can be fully exploited on current architectures, or whether they will only materialize with future systems that provide native support for low-precision arithmetic, block floating-point formats, or fine-grained precision control. This raises the broader question of how much algorithmic innovation can be decoupled from architectural co-design.

Another major challenge is automation, which relates to hardware-software co-design (\cref{sec:co-design}) and software design (\cref{sec:softwaredesign}). Manually selecting precision across levels, components, or spatial regions is not scalable. A key open question is whether the choice of precision can be automated as in~\cite{doi:10.1137/22M1522619}, based on error estimators, stability indicators, or adaptive refinement criteria, in a way that is robust, portable, and does not introduce excessive overhead. This is particularly relevant for adaptive methods, where hierarchies evolve dynamically.

Stability and robustness pose further concerns. Reducing precision locally or hierarchically alters the numerical properties of algorithms in ways that are not yet fully understood. Establishing stability guarantees, failure modes, and safe operating regimes, especially when precision varies across interacting components, remains an important open theoretical problem.

Finally, range and scaling issues become increasingly prominent as data types become narrower. Very-low-precision formats impose severe constraints on representable dynamic range, making over- and underflow central concerns. Developing analyses and algorithmic strategies that explicitly account for these range limitations, rather than assuming floating-point–like behavior, is an open and largely unexplored direction.

Together, these questions suggest that treating precision as a multilevel resource is not merely a technical optimization, but a shift in how numerical algorithms are designed, analyzed, and mapped to hardware, one whose full implications are only beginning to be understood.

\subsubsection{Multigrid}
Despite its conceptual appeal, \acrshort{CMG} with a hierarchy of precision raises several open questions that currently limit its practical adoption. A central challenge is the lack of suitable data formats and access mechanisms for extremely low precision. In particular, to show its full potential, \acrshort{CMG} requires efficient, parallel access to very small numbers (on the order of 2-4 bits per value) in \acrshort{BFP} layouts. How such values can be loaded, operated on, and stored efficiently by many threads remains an open hardware-software co-design problem.

Similarly, \acrshort{CMG} can benefit from fine bit-width gradations, such as 20, 24, or 28 bits, to match local accuracy requirements.
However, existing formats largely focus on power-of-two widths (e.g., 16, 32, 64 bits). Therefore, specialized number formats (see~\cref{sec:numberformats}) and possibly emulation 
(see~\cref{sec:emulation}) are of interest here.
Identifying which granularities are algorithmically meaningful, and whether they can be supported efficiently, is an open research question.

Another issue is the question of efficient dot products and reductions in \acrshort{BFP} formats for very long blocks, as encountered in smoothing, residual computation, and coarse-grid corrections. While short-block \acrshort{BFP} operations are increasingly supported, scalable primitives for long vectors remain lacking.
Long blocks are important in \acrshort{CMG} to keep the storage cost of exponents low in comparison to mantissas.

A further subject is the compression of operators.
Currently, \acrshort{CMG} compresses vectors, which is beneficial in the matrix-free setting.
However, certain settings, such as complex geometries, make it difficult to eliminate the need for matrix storage.
In these cases, the storage cost of matrices typically greatly dominates over the cost of vectors.
Thus, finding suitable compressed representations of system matrices and intergrid operators is essential to broaden the applicability of \acrshort{CMG} by reducing memory traffic in matrix-based solvers.
One starting point might be stencil-function representations~\cite{drzisga_surrogate_2019}.

An interesting question also arises in the context of algebraic multigrid. Algebraic multigrid consists of two phases: generation of the multigrid hierarchy and application of the multigrid method. The time and resources spent on the generation and application phases depend on the specific context. While there are several results on the effects of finite-precision errors in the application phase, we are not aware of any on the effects in the generation phase.

\subsubsection{\acrshort{AMR}}
An open research question is whether these indicators can be computed reliably in mixed or reduced precision, as a rigorous finite-precision error analysis for their evaluation is currently lacking, to our knowledge.

\section{Conclusion}
Mixed precision is no longer an optimization trick; it is a foundational method for sustaining scientific progress under energy and cost constraints. The next decade will be shaped by how well we can make low-precision trustworthy, portable, and composable---across formats, hardware, compilers, libraries, and emerging architectures. 
The core shift is conceptual: success is measured not in peak \FLOPs~but in time or energy per trusted solution, where the energy consumption is increasingly being considered the more fundamental metric for cost. 
Achieving that requires a disciplined posture---recklessly responsible computing---where aggressive approximation is paired with systematic detection, escalation, and certification. The open questions cataloged here are not just research curiosities; they define the work needed to turn mixed precision into dependable scientific infrastructure.

The central message is simple: use low precision aggressively, but pair it with mechanisms for detection, escalation to higher precision when needed, and certification of the final result. Success should ultimately be measured not by raw performance alone, but by energy per trusted solution.

\section*{Acknowledgments}
This work has benefited from Dagstuhl Seminar 26081 \say{Reduced and Mixed Precision Computing for Science and Engineering Applications}.\\ 
Research was sponsored by the Department of the Air Force Artificial Intelligence Accelerator and was accomplished under Cooperative Agreement Number FA8750-19-2-1000. The views and conclusions contained in this document are those of the authors; they should not be interpreted as representing the official policies, either expressed or implied, of the Department of the Air Force or the U.S.~Government. The U.S.~Government is authorized to reproduce and distribute reprints for Government purposes notwithstanding any copyright notation herein. \\
Bauer would like to thank the NHR-Verein e.V. (\url{www.nhr-verein.de}) for supporting this work/project within the NHR Graduate School of National High Performance Computing (NHR).\\
Carson and Ma acknowledge funding from the European Union (ERC, inEXASCALE, 101075632) and from the Charles University Research Centre program No. UNCE/24/SCI/005. \\
Fasi was supported by the Engineering and Physical Sciences Research Council [grant number UKRI2774].
D'Ambra and Vacek acknowledge funding from the European Union via the EuroHPC JU Energy oriented Center of Excellence: Fostering the European Energy Transition with Exascale (EoCoE-III), No 101144014.\\
Göddeke acknowledges support from the Stuttgart Center for Simulation Science (SimTech).\\
Graillat, J\'ez\'equel, and Mary acknowledge funding by the FPT-4 (ANR-24-CE46-7572), NumPEx Exa-MA (ANR-22-EXNU-0002), and MixHPC (ANR-23-CE46-0005) projects of the French National Agency for Research (ANR). \\
Buttari and Mary acknowledge funding by the NumPEx Exa-SofT (ANR-22-EXNU-0003) project of the French National Agency for Research (ANR). \\
Iakymchuk acknowledges funding from the European Union via the EuroHPC JU Center of Excellence in Exascale CFD (CEEC), No 101093393.\\
Kružík acknowledges funding from the European Union under the REFRESH - Research Excellence For Region Sustainability and High-tech Industries project, No CZ.10.03.01/00/22\_003/0000048 via the Operational Programme Just Transition, and under the INODIN project, No CZ.02.01.01/00/23\_020/0008487 via the Operational Programme Jan Amos Komenský.\\
Langou acknowledges funding by the NSF award \#2004850 entitled \say{\emph{Frameworks: Basic ALgebra LIbraries for Sustainable Technology with Interdisciplinary Collaboration (BALLISTIC)}}.\\
Li is supported by the U.S. Department of Energy, Office of Science, Office of Advanced Scientific Computing Research, Scientific Discovery through Advanced Computing (SciDAC) Program through the FASTMath Institute, under contract number DE-AC02-05CH11231.\\
Mikaitis acknowledges funding by the Engineering and Physical Sciences Research Council (EPSRC) grant \say{\emph{Informing Future Numerical Standards by Determining Features of Non-Standard Mathematical Hardware}}, ref. UKRI151.\\
Osei-Kuffuor acknowledges support from Lawrence Livermore National Laboratory through the Laboratory Directed Research and Development Program, project 24-ERD-035, under the auspices of the U.S. Department of Energy and under Contract DE-AC52-07NA27344.\\
Quintana acknowledges funding from projects PID2023-146569NB-C2 of MCIN/AEI/10.13039/501100011033, and CIPROM/2022/20 of the Generalitat Valenciana.\\
Vieubl\'{e} acknowledges funding from the National Natural Science Foundation of China (No. 12288201) and the Beijing Natural Science Foundation (No. IS25038).\\

\printglossary[type=\acronymtype]
\printglossary

\bibliographystyle{plain}

\bibliography{strings,references_emergingarchitectures,references,references_hierarchical_precisions,references_relatedtootherapprox,references_formats}

\end{document}